\documentclass[reqno,11pt]{amsart} 

\usepackage[top=2.0cm,bottom=2.0cm,left=3cm,right=3cm]{geometry}
\usepackage{amsthm,amsmath,amssymb}
\usepackage{mathrsfs,amsfonts,dsfont,functan,extarrows,mathtools}
\usepackage[numbers,sort&compress]{natbib}
\usepackage[colorlinks]{hyperref}

\usepackage{marginnote,multirow,makecell,booktabs}
\usepackage{xcolor}

\newtheorem{theorem}{Theorem}[section]
\newtheorem{proposition}[theorem]{Proposition}

\newtheorem{remark}[theorem]{Remark}
\newtheorem{lemma}[theorem]{Lemma}

\numberwithin{equation}{section}
\allowdisplaybreaks

\def\ls{\lesssim}
\def\gs{\gtrsim}

\def\eps{\varepsilon}

\def\p{\partial}
\def\div{{\rm div}}
\def\mix{{\rm mix}}
\def\d{\mathop{}\!\mathrm{d}} 
\def\no{\nonumber}
\def\A{\mathcal{A}}
\def\L{\mathcal{L}}
\def\P{\mathcal{P}}

\def\E{\mathcal{E}}
\def\D{\mathcal{D}}

\newcommand{\abs}[1]{\left|#1\right|}
\newcommand{\nm}[1]{\left\|#1\right\|}

\newcommand{\lx}{{L^2_x}}
\newcommand{\lxq}{{L^2_{x,q}}}
\newcommand{\hx}[1]{{H^ {#1} _x}}
\newcommand{\hxq}[1]{{H^ {#1} _x L^2_{q}}}
\newcommand{\hxqt}[1]{{H^ {#1} _{x,q}}}

\newcommand{\skpa}[2]{\left\langle #1,\, #2 \right\rangle}
\newcommand{\skpt}[2]{\langle #1,\, #2 \rangle}
\newcommand{\agl}[1]{\left\langle #1 \right\rangle}

\begin{document}

\title[Global Classical Solutions for Reactive Polymers]
{Global Existence of Classical Solutions for a Reactive Polymeric Fluid near Equilibrium} 

\author[C. Liu]{Chun Liu}
\address[Chun Liu]{\newline Department of Applied Mathematics, Illinois Institute of Technology, Chicago, IL 60616, USA}
\email{cliu124@iit.edu}

\author[Y. Wang]{Yiwei Wang}
\address[Yiwei Wang]{\newline Department of Applied Mathematics, Illinois Institute of Technology, Chicago, IL 60616, USA}
\email{ywang487@iit.edu}

\author[T.-F. Zhang]{Teng-Fei Zhang$^\dag$}
\address[Teng-Fei Zhang]{\newline School of Mathematics and Physics, China University of Geosciences, Wuhan, 430074, P. R. China
}
\email{zhangtf@cug.edu.cn}

\thanks{$^\dag$ Corresponding author.}


\begin{abstract}

In this paper, we study a new micro-macro model for a reactive polymeric fluid, which is derived recently in [Y. Wang, T.-F. Zhang, and C. Liu, \emph{J. Non-Newton. Fluid Mech.} 293 (2021), 104559, 13 pp], by using the energetic variational approach. The model couples the breaking/reforming reaction scheme of the microscopic polymers with other mechanical effects in usual viscoelastic complex fluids. 
We establish the global existence of classical solutions near the global equilibrium, in which the treatment on the chemo-mechanical coupling effect is the most crucial part. In particular, a weighted Poincar\'e inequality with a mean value is employed to overcome the difficulty that arises from the non-conservative number density distribution of each species.

\vspace*{5pt}
\noindent\textit{Keywords}: Global existence; Viscoelastic fluids; Energetic variational approach; A priori estimate; Weighted Poincar\'e inequality. 

\noindent\textit{2020 Mathematics Subject Classification}: 35A01, 35A15, 76A10, 76M30, 82D60

\end{abstract}

\maketitle

\tableofcontents


\section{Introduction} 
\label{sec:1}


\label{sub:new_micro_macro_model}



Complex fluids are fluids with complicated rheological phenomena, arising from different ``elastic'' effects, such as elasticity of deformable particles, interaction between charged ions and bulk elasticity endowed by polymer molecules. These complicated elastic effects can usually be realized by the coupling between the dynamics of macroscopic fluids and the induced elastic stress at the microscopic level. The micro-macro models have been widely used to describe the dynamics of complex fluids, including especially for polymeric fluids \cite{BCAH-87b2,DE-86b,LZ-07cms,LbL-12scm,LLZ-07cpam,Lin-12cpam}. Beyond that, for many other materials and complex fluids model like living polymers and reptation model, the mathematical theory are almost unexplored \cite{LZ-07cms}.



In this paper, we consider a reactive complex fluids model proposed firstly in \cite{LWZ-21nnfm}:
  \begin{align}\label{sys:2sp-mic-mac} \hspace*{-3.3mm}
  \begin{cases} \displaystyle
    \p_t u + u \cdot \nabla_x u + \nabla_x p = \mu \Delta_x u + \lambda \nabla_x \cdot \Big[ \int_{\mathbb{R}^3} (\nabla_q U_A \Psi_A + \nabla_q U_B \Psi_B) \otimes q \d q \Big], \\[3pt]
    \nabla_x \cdot u = 0, \\[3pt]
    \p_t \Psi_A + \nabla_x \cdot (u \Psi_A) + \nabla_q \cdot (\nabla_x u q \Psi_A) = \nabla_q \cdot (\nabla_q \Psi_A + \nabla_q U_A \Psi_A) - (k_1 \Psi_A - k_2 \Psi_B^2) , \\[3pt]
    \p_t \Psi_B + \nabla_x \cdot (u \Psi_B) + \!\!\nabla_q \cdot (\nabla_x u q \Psi_B) = \nabla_q \!\cdot\! (\nabla_q \Psi_B + \nabla_q U_B \Psi_B) + 2 (k_1 \Psi_A - k_2 \Psi_B^2).
  \end{cases}
  \end{align}
Here the unknowns $\Psi_A(t,\,x,\,q)$ and $\Psi_B(t,\,x,\,q)$ are number density distribution functions for two species $A$ and $B$ with different chain of lengths, depending on time $t\ge 0$, macroscopic position $x \in \Omega \subset \mathbb{R}^3$, and microscopic configuration variable $q \in D$ ($D$ is a bounded or unbounded domain in $\mathbb{R}^3$). $U_A(q)$ and $U_B(q)$ denote the spring potential functions.
This is a two-species model for a wormlike micellar solution (also called living polymers), in which linear chain polymers that can break and recombine continuously \cite{Cat87,VCM07,LWZ-21nnfm}. Both species of polymer particles are considered to be transported along with a incompressible fluid flow, whose evolution is described by its velocity field $u(t,\,x)$. From a viewpoint of mathematics, this model couples the incompressible Navier-Stokes equations and kinetic Fokker-Planck equations through the induced elastic stress term in \eqref{sys:2sp-mic-mac}$_1$ (the last term on the right-hand side), and the drift terms in \eqref{sys:2sp-mic-mac}$_3$-\eqref{sys:2sp-mic-mac}$_4$ (the third terms on the left-hand side).

This two-species micro-macro system \eqref{sys:2sp-mic-mac} is inspired by the macroscopic phenomenological Vasquez-Cook-McKinley (VCM) model \cite{VCM07}, 
where the longer chains can break at the middle to form shorter chains, and at the same time shorter chains can also recombine to form into longer chains. We denote the two species by $A$ and $B$, respectively, corresponding to the system \eqref{sys:2sp-mic-mac}. The breakage/reforming procedure between the two species is modeled as the chemical reaction scheme:
    $A \xrightleftharpoons[k_2]{\ k_1 \ } 2B.$
According to the \emph{law of mass action} (LMA) in chemical kinetics theory \cite{WG86}, which indicates that the rate of a reaction process is proportional to the concentrations of the reactants, the total reaction rate can be expressed as follows:
  \begin{align}\label{eq:LMA-micro}
    r= k_1 \Psi_A- k_2 \Psi_B^2,
  \end{align}
where $k_1$ and $k_2$ are the rate coefficients for breakage and reforming process respectively \cite{LWZ-21nnfm}.

It can be checked that the micro-macro model \eqref{sys:2sp-mic-mac} obeys the following entropy-entropy production law:
  \begin{align}\label{law:energy-dissipation}
    \frac{\d}{\d t} & \Big\{ \int_\Omega \tfrac{1}{2}|u|^2 \d x
      + \lambda \iint_{\Omega \times \mathbb{R}^3} \left[ \Psi_A (\ln \Psi_A + U_A -1) + \Psi_B (\ln \Psi_B + U_B -1) \right] \d q \d x
      \Big\} \no\\
    =\ & - \int_\Omega \mu |\nabla_x u |^2 \d x
         - \lambda \iint_{\Omega \times \mathbb{R}^3} \left[ \Psi_A \abs{\nabla_q (\ln \Psi_A + U_A)}^2 + \Psi_B \abs{\nabla_q (\ln \Psi_B + U_B)}^2 \right]  \d q \d x
         \no\\
       & - \lambda \iint_{\Omega \times \mathbb{R}^3} (k_1 \Psi_A - k_2 \Psi_B^2) (\ln \tfrac{\Psi_A}{\Psi_B^2} + U_A - 2U_B) \d q \d x.
  \end{align}
Indeed, the PDE system \eqref{sys:2sp-mic-mac} can be derived from this energy-dissipation law due to its variational structure, by using the \emph{energetic variational approach} (EnVarA). A detailed derivation will be given in section 2, and we also refer interested readers to \cite{LWZ-21nnfm} for details on modeling and numerical simulations. Our main concern in this present paper is the global well-posedness around the global equilibrium.

\smallskip\noindent\underline{\textbf{Brief reviews}}.
%
The study of viscoelastic fluids, as one of the most important types of complex fluids, has been an active research field that attracts more and more attentions from different areas. There are two important and famous models for dilute polymeric fluids: the FENE (Finite Extensible Nonlinear Elastic) model and the Hookean dumbbell model, up to different choices for potential functions. In the former case, the potential is chosen as $U(q) = - k \ln (1- \tfrac{|q|^2}{|b_0|^2}) $ for some constant $k>0$, in which it is commonly assumed that the polymer elongation vector belongs in a bounded ball (finite extensibility), i.e., $q \in B(0, b_0)$ with radius $b_0 >0$. In the latter Hookean case, $U(q) = \tfrac{1}{2} |q|^2$ with $q \in \mathbb{R}^3$. By considering the equations of second moments $\tau = \int_{\mathbb{R}^3} \nabla_q U \otimes q \Psi \d q$ via an approximate closure process, one can recover the Oldroyd-B model \cite{Old1950}. There are a huge of literatures in those fields, we just review here some related works of mathematical treatment in rigorous analysis, especially on well-posedness problems, and we refer interested readers to the surveys \cite{HLL-18notes,Mas-18notes} and references therein for more issues.

As for the local existence, there are many researches from different settings, see \cite{JLL-04jfa,ZZ-06arma,ELZ-04cmp,LZZ-04cpde,Ren-91sima}, for instance.
Concerning the global existence, Chemin-Masmoudi \cite{CM-01sima} proved the local and global well-posedness in critical Besov spaces. Lions-Masmoudi \cite{LM-00cam} considered the global weak solutions. We emphasize that Masmoudi proved in \cite{Mas-13invent} the global weak solutions to the FENE model by finding some new a priori estimates to justify the compactness. Hu-Lin proved recently in \cite{HL-16cpam,Hu-18jde} the global existence of weak solution for incompressible and compressible case, respectively, by assuming some small conditions on initial velocity and initial perturbated deformation tensor. As for global existence of classical solutions with small data, we refer to, for example, Lei-Liu-Zhou \cite{LLZ-08arma}, Lin-Zhang-Zhang \cite{LZZ-08cmp} and Masmoudi \cite{Mas-08cpam}. Lin-Liu-Zhang \cite{LLZ-07cpam} used the EnVarA to derive rigorously a micro-macro polymeric system and proved global existence near the global equilibrium with some assumptions on the potential $U$ (including the Hookean spring case). Corresponding to the incompressible model in \cite{LLZ-07cpam}, Jiang-Liu-Zhang \cite{JLZ-18sima} derived the compressible model by EnVarA and proved an analogous global result, where the two flow maps of both macroscopic spatial variables $x$ and microscopic configuration $q$ must be considered explicitly in the the derivation process.

The two-species micro-macro model \eqref{sys:2sp-mic-mac} is a typical example of involves kinetic-fluid/micro-macro coupling, viscoelastic effect and chemical reaction process. Precisely speaking, \eqref{sys:2sp-mic-mac} contains a micro-macro coupling between each kinetic equation of $\Psi_\alpha$ and the fluid equation of $u$, and the chemical reversible breaking/reforming reaction which satisfies the microscopic law of mass action. Furthermore, such a reversible reaction effect couples with viscoelastic effect under the presence of an extensional fluid flow. This also causes some inhomogeneous flow structures like shear banding. These above coupling effects and intrinsic features cause many difficulties for us in analyzing the well-posedness of solutions in an appropriate sense.

In addition, we mention that one can obtain from a kinetic equation the macroscopic constitutive equations by performing an approximate moment closure process. The macroscopic equations are usually referred as viscoelastic fluid model. The typical examples are the studies for classic FENE model and liquid crystal polymeric model in complex fluids theory, see \cite{HCDL-08krm,YJZ-10ccp}, for instance. Associated to the two-species micro-macro model \eqref{sys:2sp-mic-mac}, a viscoelastic system has been obtained from modelling and simulation perspectives in our another paper \cite{LWZ-21nnfm}. The method is mainly based on the maximum entropy principle, and some different approaches of closure approximations are also discussed there. This in turn raises some new topics worthy of attentions, such as the well-posedness of the viscoelastic system, and the asymptotic relation between the micro-macro model and the macroscopic viscoelastic model with respect to some scaled parameter. We will study these topics in forthcoming papers.




\noindent\underline{\textbf{Equilibrium and perturbation system}}.
Our main concern in the present paper lies in the dynamic stability analysis on the two-species micro-macro model \eqref{sys:2sp-mic-mac}. Precisely speaking, we will study the global existence of classical solution near a global equilibrium in a perturbative framework.

The \emph{global} equilibrium state of each species $\alpha= A,B$ is defined by a Maxwellian function $M_\alpha(q)$ (with a normalized coefficient $c_\alpha$ being fixed later):
  \begin{align}\label{eq:equilibrium-global}
    M_A(q) = c_A e^{-U_A}, \quad
    M_B(q) = c_B e^{-U_B}.
  \end{align}
On the other hand, the detailed balance condition for the reaction scheme requires
  \begin{align*}
    K_{\textrm{eq}} = \frac{k_2}{k_1} = \frac{M_A}{M_B^2}.
  \end{align*}
By assuming $c_A=c_B^2$ for simplicity, we have
  \begin{align}\label{eq:twice-potential}
    K_{\textrm{eq}} = \frac{k_2}{k_1} = e^{- (U_A -  2 U_B)}.
  \end{align}
As a beginning of this program, we consider here a simplified case that $K_{\rm eq} = \tfrac{k_2}{k_1} = 1$ keeps as a constant, which implies $U_A = 2 U_B$. 
Therefore, the dissipation on the reaction part in the right-hand side of \eqref{law:energy-dissipation} can be rewritten as,
  \begin{align}\label{eq:dissipation-R}
    D_\mathrm{R}
    & = \iint_{\Omega \times \mathbb{R}^3} (k_1 \Psi_A - k_2 \Psi_B^2) (\ln \tfrac{\Psi_A}{\Psi_B^2} + U_A - 2 U_B) \d q \d x \\\no
    & = \iint_{\Omega \times \mathbb{R}^3} (k_1 \Psi_A - k_2 \Psi_B^2) \ln \frac{k_1 \Psi_A}{k_2 \Psi_B^2} \d q \d x \ge 0.
  \end{align}
Furthermore, we assume $k_1 = k_2 =1$, to simplify the exposition.

One of the most important features of the two-species micro-macro model \eqref{sys:2sp-mic-mac} lies in the constraints on conservation of matter during the reversible breakage/reforming process $A \xrightleftharpoons[k_2]{\, k_1 \, } 2B$. This mass conservation is only satisfied by the quantity of total concentration $2 \Psi_A + \Psi_B$ but not any individual species, i.e.,
    $\frac{\d}{\d t} \iint_{\Omega \times \mathbb{R}^3} (2 \Psi_A + \Psi_B) \d q \d x = 0$.
Consequently, assuming initially that $\int_{\mathbb{R}^3} (2 \Psi_A + \Psi_B)|_{t=0}(x,\, q) \d q =1$, we formally get,
  \begin{align}\label{eq:mass-constraint}
    \int_{\mathbb{R}^3} (2 \Psi_A + \Psi_B)(t,\, x,\, q) \d q = 1.
  \end{align}
This will enable us to solve the constants $c_A$ and $c_B$ by combining the relation $c_A=c_B^2$.



Define the fluctuations of $\Psi_A$ and $\Psi_B$ around their global Maxwellian as follows,
  \begin{align}\label{eq:perturbation}
    \Psi_A = M_A + \sqrt{M_A} f_A, \quad \Psi_B = M_B + \sqrt{M_B} f_B.
  \end{align}
Due to the above conservation of total mass, we have:
  \begin{align*}
    \int_{\mathbb{R}^3} (2 \Psi_A + \Psi_B) \d q = 1 =  \int_{\mathbb{R}^3} (2 M_A + M_B) \d q,
  \end{align*}
which immediately implies 
  \begin{align}\label{eq:conserv-total}
    \int_{\mathbb{R}^3} (2 f_A \sqrt{M_A} + f_B \sqrt{M_B}) \d q = 0.
  \end{align}
However, due to the non-conservation for each individual perturbations $f_\alpha$, it follows,
  \begin{align}\label{eq:nonconserv-each}
    \int_{\mathbb{R}^3} f_\alpha \sqrt{M_\alpha} \d q \neq 0.
  \end{align}
This causes some difficulties from a perspective of mathematical analysis. Speaking precisely, \eqref{eq:nonconserv-each} brings a mean value term in the (weighted) Poincar\'e inequality when we try to establish a closed \emph{a-priori} estimate. This point will be specified in \S \ref{sub:strategy} later.

Inserting the perturbation scheme \eqref{eq:perturbation} into the original system \eqref{sys:2sp-mic-mac} will lead us to get a perturbative system. Direct calculations imply that,
  \begin{align}
    r = k_1 \Psi_A - k_2 \Psi_B^2
      =\ & k_1 (M_A + \sqrt{M_A} f_A) - k_1 \tfrac{M_A}{M_B^2} (M_B + \sqrt{M_B} f_B)^2 \\\no
      =\ & k_1 \sqrt{M_A} (f_A - 2 \sqrt{M_B} f_B - f_B^2),
  \end{align}
where we have used the fact $K_{\rm eq} = \tfrac{M_A}{M_B^2} = 1$. We then write
  \begin{align}\label{eq:reaction-AB}
    r_A =\ & - \tfrac{r}{\sqrt{M_A}} = - k_1 (f_A - 2 \sqrt{M_B} f_B - f_B^2), \\ \no
    r_B =\ & 2 \tfrac{r}{\sqrt{M_B}} = 2 k_1 \sqrt{M_B} (f_A - 2 \sqrt{M_B} f_B - f_B^2).
  \end{align}
Note that $k_1 = 1$ has been assumed before. On the other hand, by noticing $\nabla_q \cdot (\nabla_q \Psi_\alpha + \nabla_q U_\alpha \Psi_\alpha) = \nabla_q \cdot (M_\alpha \nabla_q (\frac{\Psi_\alpha}{M_\alpha}))$ for $\alpha = A, B$, we can define
  \begin{align}
    \mathcal{L}_\alpha f_\alpha
    =\ & \tfrac{1}{\sqrt{M_\alpha}} \nabla_q \cdot \Big[ M_\alpha \nabla_q \big( \tfrac{f_\alpha}{\sqrt{M_\alpha}} \big) \Big] 
    = \Delta_q f_\alpha + \tfrac{1}{2} \Delta_q U_\alpha f_\alpha
         - \tfrac{1}{4} |\nabla_q U_\alpha|^2 f_\alpha .
  \end{align}
Therefore, the final perturbative system for $(u,\, f_A,\, f_B)$ can be written as
  \begin{align}\label{sys:perturbative}
  \begin{cases}
    \p_t u + u \cdot \nabla_x u + \nabla_x p = \mu \Delta_x u + \lambda \div_x \int_{\mathbb{R}^3} (\nabla_q U_A \otimes q f_A \sqrt{M_A} + \nabla_q U_B \otimes q f_B \sqrt{M_B}) \d q,  \displaystyle \\[4pt]
    \div_x u = 0, \\[4pt]
    \p_t f_A + u \cdot \nabla_x f_A + \nabla_x u q \nabla_q f_A = \mathcal{L}_A f_A + r_A + \nabla_x u q \nabla_q U_A (\sqrt{M_A} + \tfrac{1}{2} f_A) ,  \\[4pt]
    \p_t f_B + u \cdot \nabla_x f_B + \nabla_x u q \nabla_q f_B = \mathcal{L}_B f_B + r_B + \nabla_x u q \nabla_q U_B (\sqrt{M_B} + \tfrac{1}{2} f_B) .
  \end{cases}
  \end{align}
Here and throughout the paper we will use the divergence symbol $\div_x$ to avoid some possible confusion between divergences on variables $x$ and $q$.

\subsection{Statement of main results} 
\label{sub:main_results}

Before we state the main results, we firstly introduce some notations.

\smallskip\noindent\underline{\textbf{Notations}}.
For simplicity, we denote by $|\cdot|_{L^2_x}$ the usual $L^2$-norm over spatial variables $x$, and by $|\cdot|_\hx{s}$ the higher-order spatial derivatives $H^s$-norm. Besides, we denote by $\nm{\cdot}_{L^2_{x,q}}$ the mixed $L^2$-norm with respect to both spatial variables $x$ and microscopic configuration variables $q$, and by $\nm{\cdot}_{H^s_{x,q}}$ the higher-order mixed derivatives. The notation $\langle \cdot,\, \cdot \rangle$ will stand for the inner product in a specific space marked by a subscript. Moreover, we will omit the subscript of space marks when it is with respect to only $x$ or both $x$ and $q$. Furthermore, the angle bracket $\langle \cdot \rangle$ denotes the integral over the configurational space, i.e., $\langle f \rangle = \int_{\mathbb{R}^3} f \d q$.

It is convenient to introduce the weighted Sobolev spaces with the weighted $L^2$-inner product defined by
  \begin{align*}
    \skpa{f}{g}_{\!_M} \triangleq \left \langle f,\, g M \right\rangle
      = \iint_{\Omega \times \mathbb{R}^3} fg \,M \d q\d x,
  \end{align*}
for any pairs $f(x,\, q),\, g(x,\, q) \in L^2_{x,\, q}$. We also use the notation $\nm{\cdot}_{M}$ to denote a norm with respect to the weighted space $L^2(M\!\d q\!\d x)$. Moreover, the subscript $M_\alpha$ in above norm is usually simplified as a single subscript $M$ when no possible confusion arises.

Let $\alpha=(\alpha_1,\, \alpha_2,\, \alpha_3) \in \mathbb{N}^3 $ be a multi-index with its length defined as $\textstyle |\alpha| = \sum_{i=1}^3 \alpha_i $. We define here the multi-derivative operator $\nabla^\alpha_x = \p_{x_1}^{\alpha_1} \p_{x_2}^{\alpha_2} \p_{x_3}^{\alpha_3}$, and we also denote $\nabla_x^k$ for $|\alpha| = k$ for simplicity. In the following texts, we will use frequently the notation $\nabla^k_l \triangleq \nabla^k_x \nabla^l_q$ to stand for the mixed derivatives over the variables $x$ and $q$.

In addition, the notation $A \ls B$ means that there exists some positive constant $C>0$ such that $A \le CB$. The notation $A \sim B$ stands for the equivalence between both sides up to a constant, i.e., some $C>0$ exists such that $C^{-1}B \le A \le CB$. Sometimes we will also omit the integral domain symbol for simplicity.



\smallskip\noindent\underline{\textbf{Main results}}.
We employ here some assumptions on the potentials $U_{\alpha}$ with $\alpha = A,\,B$, which are same as that of \cite{LLZ-07cpam}. More precisely, we assume
  \begin{align}\label{asmp-1}
    & |q| \ls (1+ |\nabla_q U_{\alpha}|), \quad \text{(sometimes we just assume $|q| \ls |\nabla_q U_{\alpha}|$ for simplicity)}, \no\\
    & \Delta_q U_{\alpha} \le C + \delta |\nabla_q U_{\alpha}|^2 \quad \text{ for } \delta < 1, \no\\ & \int_{\mathbb{R}^3} |\nabla_q U_{\alpha}|^2 M \d q \le C, \quad \int_{\mathbb{R}^3} |q|^4  M \d q \le C,
  \end{align}
and
  \begin{align}\label{asmp-2}
    & |\nabla_q^k(q \nabla_q U_{\alpha})| \ls (1+|q| |\nabla_q U_{\alpha}|), \no\\
    & \int_{\mathbb{R}^3} |\nabla_q^k(q \nabla_q U_{\alpha} \sqrt{M})|^2 \d q \le C, \\\no
    & \abs{\nabla_q^k(\Delta_q U_{\alpha} - \tfrac{1}{2} |\nabla_q U_{\alpha}|^2) } \ls (1+ |\nabla_q U_{\alpha}|^2).
  \end{align}

Define energy and energy-dissipation functionals that,
  \begin{align}\label{func:higher-deri-E}
    \E_{s}(t) & = |u|_{\hx{s}}^2 + \nm{(f_A,\,f_B)}_{\hxqt{s}}^2 + \nm{q (f_A,\,f_B)}_{\hxqt{s-1}}^2, \\[5pt]
    \label{func:higher-deri-D}
    \D_{s}(t) & = |\nabla_x u|_{\hx{s}}^2 + |\rho_A|_{\hx{s}}^2
            + \sum_{k+l \le s} \Big[ \nm{ \nabla_q \big( \tfrac{\nabla_x^k \nabla_q^l f_A}{\sqrt{M_A}} \big) }_M^2
                          + \nm{ \nabla_q \big( \tfrac{\nabla_x^k \nabla_q^l f_B}{\sqrt{M_B}} \big) }_M^2 \Big]
              \\\no & \quad
            + \sum_{k=0}^{s-1} \Big[ \nm{ q\nabla_q \big( \tfrac{\nabla_x^k f_A}{\sqrt{M_A}} \big) }_M^2
                          + \nm{ q\nabla_q \big( \tfrac{\nabla_x^k f_B}{\sqrt{M_B}} \big) }_M^2 \Big]
            + \sum_{k+l \le s} \nm{\nabla_q^l f_A - 2 \nabla_q^l f_B \sqrt{M_B}}_{\hxq{k}}^2.
  \end{align}
Here 
$\rho_\alpha = \agl{f_\alpha \sqrt{M_\alpha}} = \int_{\mathbb{R}^3} f_\alpha \sqrt{M_\alpha} \d q$ denote the macroscopic number density fluctuations.

\begin{theorem}[Global existence]\label{thm:main}
  Let $s\ge 5$. Let $(u,\, f_A,\, f_B)$ be the fluctuation near the global equilibrium $(0,\, M_A,\, M_B)$ of the two-species micro-macro model \eqref{sys:2sp-mic-mac}, with their initial data $(u_0,\, f_{A,0},\, f_{B,0})$ satisfying 
    \begin{align} \textstyle
      \Psi_{\alpha,0} = M_\alpha + \sqrt{M_\alpha} f_{\alpha,0} >0, \text{ and } \int_{\mathbb{R}^3} (2 \Psi_{A,0} + \Psi_{B,0}) \d q =1.
    \end{align}

  Then, there exists some constant $\eps$ sufficiently small, such that, if the initial fluctuation satisfies
    \begin{align}\label{eq:initial-small}
      & \E_s(0) = |u_0|^2_\hx{s} + \nm{(f_{A,0},\,f_{B,0})}^2_\hxqt{s} + \nm{(qf_{A,0},\,qf_{B,0})}^2_\hxq{s-1} \le \eps, \\\label{eq:initial-small-entropy}
      & \int_\Omega \tfrac{1}{2} |u_0|^2 \d x
        + \iint_{\Omega \times \mathbb{R}^3} \sum_{\alpha = A,\,B} (\Psi_{\alpha,0} \ln \tfrac{\Psi_{\alpha,0}}{M_{\alpha}} - \Psi_{\alpha,0} + M_{\alpha}) \d q \d x \le \eps. 
    \end{align}
  Then system \eqref{sys:2sp-mic-mac} admits a unique global classical solution $(u,\, \Psi_A,\ \Psi_B)$ with $\Psi_{\alpha} = M_\alpha + \sqrt{M_\alpha} f_{\alpha} >0$, and moreover,
    \begin{align}\label{eq:glob-energy-bdd}
      \sup_{t \in [0,+\infty)} \E_s(t) + \int_0^{+\infty} \D_s(t) \d t \le C \eps,
    \end{align}
  where the constant $C$ is independent of $\eps$.
  \end{theorem}

\begin{remark}
  In fact, as pointed out in \S \ref{sub:first_order_moments} below, the requirement of Sobolev index for first-order moments $(qf_A,\,qf_B)$ in \eqref{func:higher-deri-E}-\eqref{func:higher-deri-D} does not have to be exactly $s-1$. An index of $s' \in [4,s-1]$ suffices to close the same \emph{a-priori} estimate, see \eqref{func:higher-deri-j}.
\end{remark}

Before explaining our proof, we may explain briefly the motivation for this paper by considering the property of global equilibrium state $(0,\, M_A,\, M_B)$ of the system \eqref{sys:2sp-mic-mac}. Denote by $D(u,\, \Psi_A,\, \Psi_B)$ the dissipative part of the right-hand side in \eqref{law:energy-dissipation}, then we can infer by noticing the assumptions $k_1 = k_2 =1$ and $U_A = 2 U_B$, that
  \begin{align}
    D(u,\, \Psi_A,\, \Psi_B)
    = & \int_\Omega \mu |\nabla_x u |^2 \d x
         + \lambda \sum_{\alpha=A,B} \iint_{\Omega \times \mathbb{R}^3} \Psi_\alpha \abs{\nabla_q (\ln \tfrac{\Psi_\alpha}{M_\alpha})}^2 \d q \d x
         \\\no
       & + \lambda \iint_{\Omega \times \mathbb{R}^3} (\Psi_A - \Psi_B^2) (\ln \tfrac{\Psi_A}{\Psi_B^2}) \d q \d x.
  \end{align}
Apparently it holds $D(0,\, M_A,\, M_B) =0$. Moreover, the global equilibrium $(0,\, M_A,\, M_B)$ is a critical point of energy dissipation function $D(u,\, \Psi_A,\, \Psi_B)$.

Indeed, performing a perturbative analysis on the function $D(u,\, \Psi_A,\, \Psi_B)$ around the global equilibrium $(0,\, M_A,\, M_B)$ yields, for a small perturbative parameter $\epsilon>0$, that
  \begin{align}
    \tfrac{\d}{\d \epsilon} \big|_{\epsilon=0} D(\epsilon v,\, M_A + \epsilon \psi_A,\, M_B + \epsilon \psi_B) =0,
  \end{align}
and
  \begin{align}\label{eq:dissip-2deri}
    & \tfrac{\d^2}{\d \epsilon^2} \big|_{\epsilon=0} D(\epsilon v,\, M_A + \epsilon \psi_A,\, M_B + \epsilon \psi_B) \\[2pt] \no
    & = 2 \int_\Omega \mu |\nabla_x v |^2 \d x
         + 2 \lambda \sum_{\alpha=A,B} \iint_{\Omega \times \mathbb{R}^3} \abs{\nabla_q (\ln \tfrac{\psi_\alpha}{M_\alpha})}^2 M_\alpha \d q \d x \\\no
    & \quad + 2 \lambda \iint_{\Omega \times \mathbb{R}^3} (\tfrac{\psi_A}{\sqrt{M_A}} - 2 \psi_B)^2 \d q \d x,
  \end{align}
which keeps positive for non-trivial perturbative values.

The above analysis provides our motivation on studying the global in time existence of classical solutions near the global equilibrium state $(0,\, M_A,\, M_B)$, with small initial assumptions on perturbations, as stated in Theorem \ref{thm:main} before.

  We point out that, the right-hand side of \eqref{eq:dissip-2deri} contains three non-negative terms: fluid dissipation term, microscopic dissipation term and an additional dissipation term coming from chemical reactions. They correspond to three terms in perturbative system \eqref{sys:perturbative}: the macroscopic diffusion of Laplacian operator $\Delta_x u$ on velocity field, the kinetic diffusion of Fokker-Planck operator $\L_\alpha f_\alpha$, and the linear reaction part $(f_A - 2 \sqrt{M_B} f_B)$, respectively. This observation is very useful in constructing the dissipation functionals \eqref{func:higher-deri-D} for energy estimates.


\subsection{Strategy of the proof} 
\label{sub:strategy}


The most important part in proving the global existence Theorem \ref{thm:main} is to establish a closed \emph{a-priori} estimate. The \emph{a-priori} estimate in Proposition \ref{prop:a-priori}, together with the basic energy-dissipation law \eqref{law:energy-dissipation} (or exactly, the refined energy-dissipation law \eqref{law:refined-basic} defined below in \S \ref{sec:completion_of_proof}), and a standard bootstrap principle of continuity argument, will enable us to get the global-in-time existence result Theorem \ref{thm:main}, based on the local existence of solution under same assumptions.

Proving the \emph{a-priori} estimate requires more concerns on higher-order energy estimates for perturbations $(u,\, f_A,\, f_B)$. We firstly notice that both of energy and dissipation functionals defined in \eqref{func:higher-deri-E}-\eqref{func:higher-deri-D} contain following parts: pure spatial derivatives, spatial derivatives on first-order moment, and mixed spatial-configurational derivatives. In particular, the dissipation functionals $\D_s(t)$ in \eqref{func:higher-deri-D} is constructed with additional contributions coming from the reversible reaction scheme, which involves the linear reaction rate part $(f_A - 2 \sqrt{M_B} f_B)$, and the number density term $\rho_A$.

The whole proof is correspondingly decomposed into four steps:
  \begin{enumerate}

    \item The higher-order purely spatial derivatives for the whole perturbative system \eqref{sys:perturbative}, which contains mainly two parts of estimates, one is for macroscopic velocity field $u$ and the other for microscopic perturbative distribution functions $(f_A,\, f_B)$.

    The key point is to deal with the linear terms  arising in the micro-macro coupling: $\nabla_x \nabla_x^s u q \nabla_q U_\alpha$ and $\div_x \int_{\mathbb{R}^3} \nabla_q U_\alpha \otimes q \nabla_x^s f_\alpha \sqrt{M_\alpha} \d q$ with respect to the perturbation functions, due to the reason that a linear term is usually the worst term in the near-equilibrium-framework with small data. Fortunately, there exists a cancellation relation between the micro-macro coupling, as noticed in \cite{LLZ-07cpam}, i.e.
      \begin{align*}\textstyle
        \skpa{\nabla_x^{s+1} u q \nabla_q U_\alpha \sqrt{M_\alpha}}{\nabla_x^s f_\alpha}_{L^2_{x,q}}
        + \skpa{\div_x \int_{\mathbb{R}^3} \nabla_q U_\alpha \otimes q \nabla_x^s f_\alpha \sqrt{M_\alpha} \d q}{\nabla_x^s u}_{L^2_x} = 0.
      \end{align*}
  \end{enumerate}

However, the above step is not sufficient to close the energy estimates since some configuration moment and mixed derivatives terms are involved and hence need to be controlled (see Remark \ref{rem:esm-need}), as in the following steps:
  \begin{enumerate}
    \item[(2)] The higher-order spatial derivatives on first-order moment $qf_\alpha$: as pointed out in Remark \ref{rem:esm-need}, the highest order derivative on moment $qf_\alpha$ can be viewed only as a dissipative part, so it suffices to estimate the higher-not-highest order derivatives on first-order moment, in Sobolev space with index less than $s$.

    \item[(3)] The higher-order derivatives on mixed spatial-configurational variables $(x,q)$ is performed in a similar but more complicated way.
  \end{enumerate}

In fact, our proof in this present paper relies heavily on the Poincar\'e inequality, in a formulation of that with a Maxwellian integral weight, and with a mean value function due to the above mentioned non-conservation of perturbation, see the weighted Poincar\'e inequalities stated in Lemma \ref{lemm:Poinc-inequ-0}. This type of lower-order mean value function should be viewed only as an energy contribution, rather than a dissipation contribution equipped with one more higher-order derivatives. Notice from the simple fact that the mean value function defines actually the number density $\rho_\alpha$ of each species in a macroscopic level (see, for example, \eqref{eq:poinc-mean}), so we are required to search in the last step some potential dissipations on perturbative number density $\rho_\alpha$:
  \begin{enumerate}
    \item[(4)] The higher-order spatial derivatives on perturbative number density $\rho_\alpha$. By performing integration with respect to weighted measure $\sqrt{M_\alpha} \d q$, we can get equation \eqref{eq:number-density} of perturbative number density $\rho_\alpha$ which represents almost pure reaction contributions without any other micro-macro coupling effects, see \S \ref{sub:dissipation_density} for more details. Roughly speaking, this process enables us to separate the chemical reaction effect from a complex micro-macro coupling structure.

    At last, an additional dissipation term $(4 \agl{M_A} +1) \|\nabla_x^s \rho_A\|_{L^2_x}^2$ is obtained in \eqref{eq:dissipation-density}, due mainly to the linear reaction rate term $(f_A - 2 \sqrt{M_B} f_B)$ and the total conservation of perturbations expressed in a macroscopic formulation $2 \rho_A + \rho_B = 0$, as a direct consequence of \eqref{eq:conserv-total}.

  \end{enumerate}

In summary, combining the above four steps leads to the desired \emph{a-priori} estimate, see Proposition \ref{prop:a-priori}. As a common sense, the linear terms and lower-order terms in energy estimates are usually hard to deal with in establishing the global existence in a ``small solution'' theory. Both types of difficulties are encountered here and need more concerns. There are two types of linear terms in the perturbative system \eqref{sys:perturbative}, one of which comes from the micro-macro coupling, as we mentioned in Step (1). Our treatment is similar as that in Lin-Liu-Zhang \cite{LLZ-07cpam} and Jiang-Liu-Zhang \cite{JLZ-18sima}, and the cancellation relation will play an important role.

The other type of linear terms arise from the linear part of LMA-rate reaction contributions: $(f_A - 2 \sqrt{M_B} f_B)$ appearing in the last third term $r_\alpha$ of \eqref{sys:perturbative}. This provides a contribution to dissipation in two different scales. At a microscopic level, this exhibits an explicit formulation of a positive quadratic dissipation term $\nm{\nabla_x^k \nabla_q^l f_A - 2 \sqrt{M_B} \nabla_x^k \nabla_q^l f_B }_{L^2_x}^2$, as appeared in the last term in the definition of dissipation functional \eqref{func:higher-deri-D}. At a macroscopic level, there is an implicit formulation through the dissipation of number density $(4 \agl{M_A} +1) \|\nabla_x^s \rho_A\|_{L^2_x}^2$, as we have stated in above Step (4). Note that, this linear reaction dissipation contribution seems more important on a macroscopic level rather than that on a microscopic level. Besides, one can find in \cite{LWZ-21nnfm} some similar treatment on the macroscopic number density, from a viewpoint of modeling.

  At the end, we point out that the treatment for number density in step (4) can also be understood in a terminology of kinetic theory, especially of the hydrodynamic limit theory for the Boltzmann equation, see for instance Guo \cite{Guo-02cpam,Guo-03invent,Guo-06cpam} and Liu-Yu-Yang \cite{LY-04cmp,LYY-04physd}. Because the number density lies actually in the kernel of Fokker-Planck operator which has only one dimension spanned by constant variable (up to the weighted measure), the above integral process can be viewed as the macroscopic projection on kernel space of Fokker-Planck operator. On the other hand, the Poincar\'e inequality applied to the pure microscopic part will show a clear formulation without a mean value. There is a close relation between our treatment here and the so-called macro-micro decomposition in hydrodynamic limit theory. However, the non-conservative kinematics \eqref{eq:kinematic} below and the LMA-rate reaction appeared in the basic energy-dissipation law \eqref{law:energy-dissipation} make an effort in energy estimates, both of which are very different from the classical kinetic theory.



\subsection{Organization of the paper} 
\label{sub:organization_of_the_paper}


The rest of this paper is as follows: a formal derivation of this two-species micro-macro model \eqref{sys:2sp-mic-mac} will be given, by using the EnVarA in the following section \S \ref{sec:derivations}, containing derivations for the mechanical and chemical reaction parts.


In \S \ref{sec:lemmas}, some useful lemmas are presented, involving the structure of the Fokker-Planck operator and the weighted Poincar\'e inequalities. \S \ref{sec:a_priori_estimates} is devoted to closing the \emph{a priori} estimate, see Proposition \ref{prop:a-priori}, for that we need to derive closed estimates for fluctuations, including purely spatial higher-order derivatives on fluctuations themselves and fluctuation moments, mixed spatial-configurational derivatives, and furthermore, additional dissipation on number densities at the macroscopic level. In the last section \S \ref{sec:completion_of_proof}, we complete the proof of global existence result (Theorem \ref{thm:main}), by combining Proposition \ref{prop:a-priori}, the basic energy-dissipation law, and a standard bootstrap principle of continuity argument together.

\section{A Formal Energetic Variational Derivation} 
\label{sec:derivations}

In this section, we provide a formal derivation of our model \eqref{sys:2sp-mic-mac} by the \emph{Energetic Variational Approach} (EnVarA), in the spirit of \cite{LWZ-21nnfm}. We present it here in a self-contained way, for the sake of readability of this paper. The framework of EnVarA, which was developed from seminal work of Rayleigh \cite{Ray1871} and Onsager \cite{Ons1931-1,Ons1931-2}, has been a powerful tool to deal with the couplings and competitions between different mechanisms in different scales, such the competition between the microscopic interactions and the macroscopic dynamics in complex fluids (see the surveys \cite{Lc-09notes,GKL-18notes}, for example).

Starting from the energy-dissipation law and the kinematic relations, EnVarA provides a unique, well-defined way to derive the dynamics of the systems (which is usually represented as a coupled system of nonlinear partial differential equations). The main ingredients include the \emph{Least Action Principle} (LAP) and the \emph{Maximum Dissipation Principle} (MDP), which derive the conservative force and the dissipative force respectively. The force balance condition will lead to the final PDE system. The EnVarA has been successfully applied to model many systems, especially those in complex fluids, such as liquid crystals, polymeric fluids, phase field and ion channels \cite{Lc-09notes,GKL-18notes}. Since these models are derived based on the basic thermodynamic laws and hence are thermodynamically consistent. Besides, this approach can also be used to study problems with boundary, especially for dynamical boundary conditions problems for the Cahn-Hilliard equation \cite{LW-19arma,KLLM-20m2an}.

Recently, EnVarA has been applied in the study of reaction-diffusion systems driven by the \emph{law of mass action} (LMA), see Wang-Liu-Liu-Eisenberg \cite{WLLE-20pre}. The authors applied a generalized notion of EnVarA to a non-equilibrium reaction-diffusion system, by finding a way to couple successfully the chemical reaction with other effects. It is worthy to point out that their treatment is different from the \emph{linear response theory} in which the dissipation rate function of total energy is of a quadratic formulation, as a common assumption in non-equilibrium thermodynamic theory \cite{Ons1931-1,Ons1931-2}. The authors showed in \cite{WLLE-20pre} that the dynamics of system is determined by the choice of the dissipation.

In this paper, we consider a micro-macro model for such wormlike micellar solutions involving the chemical reaction of breakage/reforming scheme. The kinematic relations for the two species can be written as follows:
  \begin{align}\label{eq:kinematic}
  \begin{cases}
    \p_t \Psi_A + \nabla_x \cdot (u_A \Psi_A) + \nabla_q \cdot (V_A \Psi_A) = -r , \\
    \p_t \Psi_B + \nabla_x \cdot (u_B \Psi_B) + \nabla_q \cdot (V_B \Psi_B) = 2 r.
  \end{cases}
  \end{align}
where $u_\alpha$ and $V_\alpha$ are the corresponding macroscopic and microscopic velocities in spatial and configurational spaces, respectively. Since we consider the case the molecule orientations $q$ transporting along with the environmental fluid velocity $u(t,\,x)$, so $u_A= u_B = u$. Besides, the fluid velocity field $u$ is considered as an incompressible flow. Note that the left-hand side of each equation represents a formulation of Smoluchowski equation \cite{BCAH-87b2,DE-86b}, while the right-hand side represents the chemical reaction of breakage/reforming process between the two species, and moreover, $r=k_1 \Psi_A- k_2 \Psi_B^2$ by recalling \eqref{eq:LMA-micro}.

Based on the kinematics \eqref{eq:kinematic}, the two-species micro-macro model can be derived from the following basic energy-dissipation law:
  \begin{align}\label{law:general-diffu}
    \frac{\d}{\d t} & \Big\{ \int_\Omega \tfrac{1}{2}|u|^2 \d x
    + \lambda \sum_{\alpha=A,B} \iint_{\Omega \times \mathbb{R}^3} \Psi_\alpha (\ln \Psi_\alpha + U_\alpha -1) \d q \d x
    \Big\} \no\\
    =\ & - \int_\Omega \mu |\nabla_x u |^2 \d x
         - \lambda \sum_{\alpha=A,B} \iint_{\Omega \times \mathbb{R}^3} \Psi_\alpha \abs{V_\alpha - \nabla_x u q}^2 \d q \d x
         \no\\
       & - \lambda \iint_{\Omega \times \mathbb{R}^3} (k_1 \Psi_A - k_2 \Psi_B^2) (\ln \tfrac{\Psi_A}{\Psi_B^2} + U_A - 2 U_B) \d q \d x,
  \end{align}
by employing the EnVarA \cite{Lc-09notes,GKL-18notes}. The second dissipative term in the right-hand side stands for the relative friction of microscopic polymer particles to the environmental macroscopic flow, in which the expression $\nabla_x u q$ comes from the Cauchy-Born rule $q=FQ$ with $F$ being the deformation tensor and $Q$ being the initial configuration in Lagrangian coordinates.

\subsection{The mechanical part} 
\label{sub:the_mechnaical_part}


Noticing the variational structure of above energy-dissipation law, it is convenient to employ the EnVarA, which has the advantage of describing mathematically the competitions between the kinetic energy and free energy (including both entropy and internal energy, or say, elastic energy here). We firstly point out that both Lagrangian and Eulerian coordinates are used in this approach, where the flow map for macroscopic position $x=x(t,\ X)$, associated to the Lagrangian coordinate of initial position $X$ and a given velocity field $u(t,\,x)$, is defined as
  \begin{align}
  \begin{cases}
    \tfrac{\d}{\d t} x =u(t,\, x(t,\ X)), \\
    x(t,\ X)|_{t=0}= X.
  \end{cases}
  \end{align}
Similarly, associated to its Lagrangian coordinate $Q$ mentioned earlier, the flow map $q=q(t,\,X,\, Q)$ in the configuration space, generated by a given microscopic velocity $V=V(t,\,x,\, q)$, is defined through the formula $\tfrac{\d}{\d t} q =V(t,\, x(t,\ X),\, q(t,\,X,\, Q))$. We also introduce the deformation tensor describing the elasticity: $F(t,\, X)= \tfrac{\p x}{\p X}$ (in Lagrangian coordinate) and $F(t,\,x(t,\, X)) = F(t,\, X)$ (in Eulerian coordinate).

To apply the EnVarA, we firstly write the energy functional in Lagrangian coordinates, corresponding to its kinetic energy and free energy part, that,
  \begin{align}
    \mathcal{K} + \mathcal{F}
    & = \int_{\Omega^t} \tfrac{1}{2} |\dot x|^2 \d x
        + \lambda \sum_{\alpha=A,B} \iint_{\Omega^t \times \mathbb{R}^3} \Psi_\alpha (\ln \Psi_\alpha + U_\alpha(q) - 1) \d q \d x \no\\
    & = \int_{\Omega^0} \tfrac{1}{2} |\dot x|^2 \d X
        + \lambda \sum_{\alpha=A,B} \iint_{\Omega^0 \times \mathbb{R}^3} \Psi_{\alpha,0}(X,Q) \left( \ln \frac{\Psi_{\alpha,0}(X,Q)}{\det \tfrac{\p q}{\p Q}} + U_\alpha(q) - 1 \right) \d Q \d X,
  \end{align}
where we have used the fact $\det F=1$ due to the incompressibility of fluid flow. This formulation will allow us to define the least action functional $\mathcal{A} = \int_0^T (\mathcal{K} - \mathcal{F}) \d t$.

Note that the notion of ``separation of scales'' plays a crucial rule in the derivation process. In fact, the microscopic configuration variable $q$ and $V$ should be viewed independent from macroscopic variable $x$, which means roughly that the fluid flow is like a ``static'' background. By the LAP, performing the variation with respect to the configuration variable $q$, the variation on $\int_0^T \mathcal{K} \d t$ will vanish and we get
  \begin{align}
    \delta_q \int_0^T \mathcal{F} \d t
    =\ & \iiint \Bigg\{ \Psi_{\alpha,0} \frac{\det \tfrac{\p q}{\p Q}}{\Psi_{\alpha,0}} \cdot
                \Big[ -\tfrac{\Psi_{\alpha,0}}{(\det \tfrac{\p q}{\p Q})^2} \cdot
                        \det \tfrac{\p q}{\p Q} \cdot \text{tr} \left( \tfrac{\p Q}{\p q} \cdot \tfrac{\p \delta q}{\p Q} \right)
                \Big] \no\\[4pt] & \hspace*{1cm}
       + \Psi_{\alpha,0} \nabla_q U_{\alpha} \cdot \delta q \Bigg\} \d Q \d X \d t \no\\
    =\ & \int_0^T \iint_{\Omega^0 \times \mathbb{R}^3} \left[ - \Psi_{\alpha,0} (\nabla_q \cdot \delta q) + \Psi_{\alpha,0} \nabla_q U_{\alpha} \cdot \delta q \right]  \d Q \d X \d t \\\no
    =\ & \int_0^T \iint_{\Omega^t \times \mathbb{R}^3} \left( \nabla_q \Psi_{\alpha} \cdot \delta q + \Psi_{\alpha} \nabla_q U_{\alpha} \cdot \delta q \right)  \d q \d x \d t \\\no
    =\ & \skpa{ \nabla_q \Psi_{\alpha} + \nabla_q U_{\alpha} \Psi_{\alpha} }{\delta q}_{L^2_{t,x,q}}.
  \end{align}
On the other hand, the MDP, meaning to perform variation on dissipation part with respect to the microscopic ``velocity'' $V = \dot q$, yields,
  \begin{align}
    \delta_{V} \mathcal{D}_\mathrm{q}
    = \iint \Psi_{\alpha} (V_{\alpha} - \nabla_x u q) \cdot \delta V_{\alpha} \d q \d x
    = \skpa{\Psi_{\alpha} (V_{\alpha} - \nabla_x u q)}{\delta V_{\alpha}}_{L^2_{x,q}}.
  \end{align}
So we can obtain the force balance
  \begin{align}
    \Psi_{\alpha} (V_{\alpha} - \nabla_x u q) = - (\nabla_q \Psi_{\alpha} + \nabla_q U_{\alpha} \Psi_{\alpha}),
  \end{align}
i.e.,
  \begin{align}\label{eq:force-balance-q}
    V_{\alpha} - \nabla_x u q = - \nabla_q (\ln \Psi_{\alpha} + U_{\alpha}).
  \end{align}

With this relation in hand, we are able to consider the energetic variation on macroscopic level to get the momentum equation. The LAP implies, by taking variation over $x$, that
  \begin{align}
    \delta_x \int \mathcal{K} \d t = \iint \dot{x} \cdot \tfrac{\d}{\d t} \delta x \d X \d t
    = - \iint \tfrac{\d }{\d t}u \cdot \delta x \d x \d t
    = \skpt{-(\p_t u + u \cdot \nabla_x u)}{\delta x}_{L^2_{t,x}},
  \end{align}
while $\delta_x \int \mathcal{F} \d t$ vanishes because of the notion ``separation of scales''. Meanwhile, the MDP yields, by taking variation over $u = \dot{x}$, that
  \begin{align}\label{eq:D-vari-u}
    \delta_{u} \mathcal{D}_\mathrm{u}
    =\ & \int \mu \nabla_x u \nabla_x (\delta u) \d x + \lambda \sum_{\alpha=A,B} \iint \Psi_{\alpha} (V_{\alpha} - \nabla_x u q)\cdot (-\nabla_x \delta u) q \d q \d x
    \no\\
    =\ & \skpa{-\mu \Delta_x u}{\delta u}_{L^2_x} + \textstyle \lambda \sum_{\alpha} \skpa{\nabla_x \cdot \int \Psi_{\alpha} (V_{\alpha} - \nabla_x u q) \otimes q \d q}{\delta u}_{L^2_x}
    \no\\
    =\ & \skpa{-\mu \Delta_x u}{\delta u}_{L^2_x} - \textstyle \lambda \sum_{\alpha} \skpa{\nabla_x \cdot \int (\nabla_q \Psi_{\alpha} + \nabla_q U_{\alpha} \Psi_{\alpha}) \otimes q \d q}{\delta u}_{L^2_x}
    \no\\
    =\ & \skpa{-\mu \Delta_x u}{\delta u}_{L^2_x} - \textstyle \lambda \sum_{\alpha} \skpa{\nabla_x \cdot \left( \int \nabla_q U_{\alpha} \otimes q \Psi_{\alpha} \d q \right) - \nabla_x n_{\alpha}}{\delta u}_{L^2_x},
  \end{align}
where we have used the above microscopic force balance \eqref{eq:force-balance-q} in the penultimate line. Therefore, we get the momentum equation that,
  \begin{align}\label{eq:momentum}
    \p_t u + u \cdot \nabla_x u + \nabla_x p = \mu \Delta_x u + \lambda \nabla_x \cdot \Big[ \int_{\mathbb{R}^3} (\nabla_q U_A \Psi_A + \nabla_q U_B \Psi_B) \otimes q \d q \Big],
  \end{align}
with the pressure $p$ being the Lagrangian multiplier.

\begin{remark}
We point out that, the above induced elastic stress (Kramer's stress), namely, $\int_{\mathbb{R}^3} \nabla_q U_{\alpha} \otimes q\, \Psi_{\alpha} \d q$, is able to be derived either by the MDP (see \eqref{eq:D-vari-u}), or by the LAP by virtue of the Cauchy-Born rule, as did in \cite{LLZ-07cpam,JLZ-18sima}. This means that the elastic stress can be viewed as dissipative or conservative. The same idea can be found in the study of liquid crystal, for example, \cite{WXL-13arma}.

\end{remark}

\subsection{The chemical reaction part} 
\label{sub:the_chemical_reaction_part}


In order to take the special reaction dissipation into account, here we need a generalized EnVarA, by employing the \emph{reaction trajectory} $R$ defined through
  \begin{align}
    \tfrac{\d}{\d t} R =r,
  \end{align}
here $r$ is the LMA-rate defined in \eqref{eq:LMA-micro}. We refer readers to \cite{OP-74arma,WLLE-20pre} for more detailed discussion. Note that the reaction trajectory $R$ in the chemical reaction system is analogy to the flow map $x(X, t)$ in the mechanical system. Due to the ``separation of scales'' again, the chemical reaction of reversible breakage/reforming process $A \xrightleftharpoons[k_2]{\, k_1 \, } 2B$ can be expressed as,
  \begin{align}
    \Psi_A = - R + \Psi_{A,0}, \quad \Psi_B = 2R + \Psi_{B,0}.
  \end{align}
Actually, this relation may be regarded as the kinematics for the above breakage/reforming reaction process.

The new state variable of reaction trajectory enables us to rewrite the free energy $\mathcal{F}$ as,
  \begin{align}\label{eq:energy-R}
    \mathcal{F} (R) = \mathcal{F} (\Psi_A(R),\,\Psi_B(R)).
  \end{align}
Meanwhile, we also write the reaction dissipation part \eqref{eq:dissipation-R} as $\mathcal{D}_\mathrm{R} = \mathcal{D}_\mathrm{R} (R,\, \dot{R})$. So the energy-dissipation law for a pure chemical reaction process can be rewritten as
  \begin{align}\label{eq:EnDis-law-R}
    \frac{\d}{\d t} \mathcal{F} (R) = - \mathcal{D}_\mathrm{R} (R,\, \dot{R}).
  \end{align}

We assume that the nonnegative reaction dissipation $\mathcal{D}_\mathrm{R}$ takes the form of
  \begin{align}
    \mathcal{D}_\mathrm{R} (R,\, \dot{R}) = \skpa{\mathcal{G}_\mathrm{R}(R,\, \dot{R})}{\dot{R}},
  \end{align}
this, combining with the fact $\tfrac{\d}{\d t} \mathcal{F} (R) = \skpa{\tfrac{\delta \mathcal{F}}{\delta R}}{\dot{R}}$, yields that,
  \begin{align}\label{eq:grad-flow-gnrl}
    \mathcal{G}_\mathrm{R} (R, \dot{R}) = - \frac{\delta \mathcal{F}}{\delta R}.
  \end{align}
We refer this as a general gradient flow. The exact expression of LMA-rate $r$ will be revisited by choosing
  \begin{align}
    \mathcal{D}_\mathrm{R} (R,\, \dot{R}) = \dot{R} \ln \left( \tfrac{\dot{R}}{k_2 \Psi_B^2} + 1 \right),
  \end{align}
Indeed, direct calculations imply that
  \begin{align}
    \frac{\delta \mathcal{F}}{\delta R}
    = \sum_{\alpha}  \frac{\delta \mathcal{F}}{\delta \Psi_{\alpha}} \cdot \frac{\p \Psi_{\alpha}}{\p R}
    & = - (\ln \Psi_A + U_A) + 2 (\ln \Psi_B + U_B) \no\\
    & = - (\ln \tfrac{\Psi_A}{\Psi_B^2} + U_A - 2 U_B)
    = - \ln \tfrac{k_1 \Psi_A}{k_2 \Psi_B^2},
  \end{align}
where we have used the equilibrium relation \eqref{eq:twice-potential}. As a result, the above general gradient flow \eqref{eq:grad-flow-gnrl} implies,
  \begin{align}
    r = \dot{R} = k_1 \Psi_A - k_2 \Psi_B^2,
  \end{align}
which is exactly the same formulation as \eqref{eq:LMA-micro}.

Therefore, inserting the above LMA-rate formulation, and the microscopic force balance relation \eqref{eq:force-balance-q} into the kinematics \eqref{eq:kinematic} enables us to get, for the number density distribution functions $(\Psi_A,\, \Psi_B)(t,\,x,\,q)$, that
  \begin{align}\label{sys:Micro-macro-2}
  \begin{cases}
    \p_t \Psi_A + \nabla_x \cdot (u \Psi_A) + \nabla_q \cdot (\nabla_x u q \Psi_A) = \nabla_q \cdot (\nabla_q \Psi_A + \nabla_q U_A \Psi_A) - (k_1 \Psi_A - k_2 \Psi_B^2) , \\
    \p_t \Psi_B + \nabla_x \cdot (u \Psi_B) + \nabla_q \cdot (\nabla_x u q \Psi_B) = \nabla_q \cdot (\nabla_q \Psi_B + \nabla_q U_B \Psi_B) + 2 (k_1 \Psi_A - k_2 \Psi_B^2).
  \end{cases}
  \end{align}

Finally, we have derived the whole two-species micro-macro model \eqref{sys:2sp-mic-mac} for wormlike micelles by combining the momentum equation \eqref{eq:momentum} and above equations for the number density distribution functions. We refer the readers to \cite{LWZ-21nnfm} for a more detailed derivation on model \eqref{sys:2sp-mic-mac} and its moment closure model.


\subsection{Comments on the derived model} 
\label{sub:comments}


We make some comments here, on some general case of equilibriums and on related models including especially the VCM model.

\smallskip\noindent\underline{\textbf{General cases}}.
In a general setting, the chemical equilibrium ``constant'' $K_{\rm eq}$ might be possibly dependent upon the variables $q$ (and moreover, upon the variables $x$). This requires more complex treatment both from analysis and simulation viewpoints.

Indeed, taking the classic Hookean elastic springs for example, namely, $U_\alpha= \tfrac{1}{2} H_\alpha |q|^2$ for $\alpha = A,B$. As mentioned above, the simplest case $K_{\rm eq}=1$ combined with equation \eqref{eq:twice-potential} yields $U_A = 2U_B$, and hence $H_A = 2 H_B$. This stands for the parallel connection when two shorter polymeric particles $B$ recombine themselves to form one longer $A$. Another physical case interested is the series connection, i.e., $H_B = 2 H_A$, as studied in VCM model \cite{VCM07}. For that, we have to deal with a non-trivial dependence of $K_{\rm eq}= K_{\rm eq} (q)$ on variables $q$, which needs more treatment for the LMA reaction rate term \eqref{eq:LMA-micro}, $r= k_1 \Psi_A- k_2 \Psi_B^2$, with non-trivial rate coefficients $k_1(q)$ and $k_2(q)$.

Furthermore, it is possible that rate coefficients $k_1$, $k_2$ (and the equilibrium ``constants'' $K_{\rm eq}$) depend not only on variables $q$ but also on variables $x$, when the local equilibrium states can also be considered, see \cite{LWZ-21nnfm} for details.

We also mention that, a compatibility condition had been employed, requiring that the chemical reaction equilibrium state is consistent with the equilibrium of each species themselves, see \eqref{eq:equilibrium-global}-\eqref{eq:twice-potential}. For now, more general inconsistent cases still remain open.

\smallskip\noindent\underline{\textbf{Related models}}.
There exist various models describing wormlike micellar solutions.
One famous model is Cates' living polymer theory. In \cite{Cat87}, Cates studied reptation dynamics of the reversible breaking/reforming reaction process of the micellar chains, in which the polymers with chain of length $L$ is assumed to be continuous and can break anywhere into polymers of shorter lengths. Based on a discrete version of Cates' model, Vasquez-Cook-McKinley (VCM) \cite{VCM07} proposed a simplified model where only two different species are incorporated to describe the flow behavior of wormlike micellar solutions coupled of a viscoelastic fluid rheology. Later, German-Cook-Beris \cite{GCB-13thermo} revisited VCM model via a general bracket approach from a viewpoint of non-equilibrium thermodynamics.

The two-species micro-macro model \eqref{sys:2sp-mic-mac} in this paper is mainly inspired by the VCM model \cite{VCM07}.
The main difference between our model and VCM model lies on the assumption for the microscopic reaction mechanism, which enables us to derive the model with a clear variational structure. To our best knowledge, there is no any statement about this energy/entropy structure in VCM model. The VCM model assumes a convolution form of reaction rate for the reforming process $\Psi_B * \Psi_B$, which leads to the kinematics involving LMA-rate at a macroscopic level:
  \begin{align}
  \begin{cases}
    \p_t n_A + \nabla_x \cdot (u_A n_A) = - k_1 n_A + k_2 n_B^2 , \\
    \p_t n_B + \nabla_x \cdot (u_B n_B) = 2 (k_1 n_A - k_2 n_B^2),
  \end{cases}
  \end{align}
where $n_\alpha = \int_{\mathbb{R}^3} \Psi_\alpha \d q$, ($\alpha = A,B$), is the number density for each species. So the VCM model is indeed a macroscopic model, due to the macroscopic LMA formulation $\tilde r = k_1 n_A - k_2 n_B^2$ on the right-hand side. In contrast of the VCM model, the LMA-rate in our model, $r= k_1 \Psi_A- k_2 \Psi_B^2$ in \eqref{eq:LMA-micro}, is assumed in a microscopic level.

On the other hand, we also mention there are some other related models concerning similar chemical reaction process, such as coagulation/fragmentation, merging/splitting model in different setting, see \cite{GST-13siap,DLP-17jns} and references therein.

\section{Preliminaries and Lemmas} 
\label{sec:lemmas}


\subsection{The kernel structure of Fokker-Planck operator} 
\label{sub:kernel_of_FP}


Let us begin with the linear Fokker-Planck operator $\mathcal{L}_\alpha f_\alpha = \frac{1}{\sqrt{M_\alpha}} \nabla_q \cdot \left[ M_\alpha \nabla_q \left( \frac{f_\alpha}{\sqrt{M_\alpha}} \right) \right]$, which naturally defines the inner product, (the subscripts $A$, $B$ are omitted here)
  \begin{align}
    \skpa{\mathcal{L} f}{f}_{L^2_q}
    & = \skpa{\tfrac{1}{\sqrt M} \nabla_q \cdot \big[ M \nabla_q ( \tfrac{f}{\sqrt M} ) \big]}{f}_{L^2_q} \\\no
    & = \skpa{M \nabla_q ( \tfrac{f}{\sqrt M} )}{\nabla_q ( \tfrac{f}{\sqrt M} ) }_{L^2_q}
    = - \int \abs{\nabla_q (\tfrac{f}{\sqrt M}) }^2 M \d q .
  \end{align}
This fact can be recast, by denoting $f = g \sqrt M$ and employing a new weighted inner product $\skpa{F}{G}_{L^2_M} = \skpa{F}{GM}_{L^2_q}$, as follows,
  \begin{align}
    \skpa{\mathcal{A} g}{g}_{L^2_M} = \skpa{\nabla_q g}{\nabla_q g}_{L^2_M} = - \int |\nabla_q g|^2 M \d q = - \abs{\nabla_q g}^2_{L^2_M},
  \end{align}
where the operator $\mathcal{A}$ is defined as $\mathcal{A} g = \frac{1}{M} \nabla_q \cdot (M \nabla_q g)$, and we have used $L^2_q$ and $L^2_M$ as a shorthand for $L^2(\d q)$ and $L^2(M \!\d q)$, respectively.

As a well-known property of Fokker-Planck operator, we state the following lemma (see \cite{DL-97arma,DFT-10cmp} for instance):
\begin{lemma}
Note that we have
\begin{enumerate}
  \item Both of $\mathcal{A}$ and $\mathcal{L}$ are self-adjoint operators with respect to their scalar products on $L^2(M \!\d q)$ and $L^2(\d q)$, respectively;
  \item Kernel of $\mathcal{A}$ is $1$-dimensional that ${\rm Ker}\mathcal{A} = {\rm span} \{ 1 \}$ on $L^2(M \!\d q)$.
\end{enumerate}
\end{lemma}

Next we define the macroscopic quantities of number density and first-order moment as
  \begin{align}
    \rho^g(t,x) = \int g(t,x,q) M \d q = \agl{g}_M, \quad
    j^g(t,x) = \int q g(t,x,q) M \d q =\agl{qg}_M.
  \end{align}


\smallskip\noindent\underline{\textbf{The coercivity estimates}}:

Define the projections on the kernel of Fokker-Planck operator $\P_0 g = \rho^g$. The fact
  \begin{align}
    \mathcal{A} g = \A (\P_0 g + \P_0^\perp g) = \A \P_0^\perp g,
  \end{align}
together with the self-adjoint property, yields
  \begin{align}
    \skpa{\mathcal{A} \P_0^\perp g}{\P_0^\perp g}_{M}
    = - \abs{\nabla_q \P_0^\perp g}^2_{L^2_M}.
  \end{align}

On the other hand, recall the weighted Poincar\'e inequality with a mean value, that
  \begin{align}\label{eq:poinc-mean}
    \int |\nabla_q h|^2 M \d q \ge \lambda_0 \int |h - \overline{(h)}_M|^2 M \d q,
  \end{align}
with the mean value defined as $\overline{(h)}_M = \frac{\int h M \d q}{\int M \d q}$. Then we can derive, for a kernel orthogonal function $h = \P_0^\perp g$, the following \emph{partial coercivity estimate} that
  \begin{align}\label{eq:partial-coer}
    \skpa{-\mathcal{A} g}{g}_{M} = \iint |\nabla_q \P_0^\perp g|^2 M \d q \d x \ge \lambda_0 \iint \abs{\P_0^\perp g}^2 M \d q \d x.
  \end{align}
Due to the commutative property of $\nabla_x$ and $\A$, we have the similar higher-order coercivity estimate:
  \begin{align}
    \skpa{-\mathcal{A} \nabla_x^s g}{\nabla_x^s g}_{M} = \iint |\nabla_q \P_0^\perp (\nabla_x^k g)|^2 M \d q \d x \ge \lambda_0 \iint \abs{\P_0^\perp (\nabla_x^k g)}^2 M \d q \d x.
  \end{align}
We mention here that the coercivity estimates are of crucial importance in establishing additional dissipation effect for perturbative number density, see \S \ref{sub:dissipation_density} for more details.

\subsection{The weighted Poincar\'e inequality with mean value} 
\label{sub:poincar'e_inequality_with_mean_value}


We state here the following lemmas describing the weighted Poincar\'e inequality with a mean value. The key part of its proof is to transfer the considered space $L^2(\d q \!\d x)$ into a weighted one $L^2(M\!\d q \!\d x)$ with respect to a Maxwellian weight $M=M(q)$. The proofs are similar as that for Lemmas 3.2--3.3 in \cite{LLZ-07cpam} (or Lemma 1.6 in \cite{JLZ-18sima}), except that an additional mean value term appeared in the present case. We thus omit the proofs here.

\begin{lemma}[Weighted Poincar\'e Inequality]\label{lemm:Poinc-inequ-0}
  We have the weighted Poincar\'e inequalities with mean value over both variables $x,\,q$, that
  \begin{align}
    \nm{f}_{\lxq}^2 & 
      \ls \iint \abs{ \nabla_q \left( \tfrac{f}{\sqrt M} \right) }^2 M \d q \d x + \abs{\rho}_{\lx}^2, \\
    \nm{\nabla_q U f}_{\lxq}^2 & 
      \ls \iint \abs{ \nabla_q \left( \tfrac{f}{\sqrt M} \right) }^2 M \d q \d x + \abs{\rho}_{\lx}^2, \\
    \nm{q\nabla_q U f}_{\lxq}^2 & 
      \ls \iint \abs{ \agl{q} \nabla_q \left( \tfrac{f}{\sqrt M} \right) }^2 M \d q \d x + \abs{\rho}_{\lx}^2,
  \end{align}
where we have used the notation $\agl{q} = (1+ |q|^2)^{1/2}$.

\begin{remark}
  We make two remarks here:
  \begin{enumerate}

    \item The proof relies on weighed Poincar\'e inequalities (such as \eqref{eq:poinc-mean}), in which we will use the fact for the mean value function, that
      \begin{align*}
        \abs{\overline{(g)}_M}_{\lx}^2
        = \abs{\frac{\int g M \d q}{\int M \d q}}_{\lx}^2
        = \agl{M}^{-1} \abs{\rho^g}_{\lx}^2.
      \end{align*}

    \item Note that these inequalities can be directly extended to the case of pure spatially higher-order derivatives, by replacing $f$, $g$ and $\rho$ by $\nabla_x^k f$, $\nabla_x^k g$ and $\nabla_x^k \rho$, respectively.
  \end{enumerate}
\end{remark}

\end{lemma}

We next state the following lemma for similar inequalities involving higher-order mixed derivatives:
\begin{lemma}\label{lemm:Poinc-inequ-mix}
  We have, for higher-order mixed derivatives with $l \ge 1$, that
  \begin{align}
    \nm{\nabla_q U \nabla_x^k \nabla_q^l f}_{\lxq}^2
      & \ls \sum_{m=0}^{l} \iint \abs{ \nabla_q \left( \tfrac{\nabla_x^k \nabla_q^m f}{\sqrt M} \right) }^2 M \d q \d x + \abs{\nabla_x^k \rho}_{\lx}^2,  \\[5pt]
    \nm{\nabla_x^k \nabla_q^l f}_{\lxq}^2
      & \ls \sum_{m=0}^{l-1} \iint \abs{ \nabla_q \left( \tfrac{\nabla_x^k \nabla_q^m f}{\sqrt M} \right) }^2 M \d q \d x + \abs{\nabla_x^k \rho}_{\lx}^2.
  \end{align}
\end{lemma}


\section{The A-Priori Estimates} 
\label{sec:a_priori_estimates}



\subsection{The a-priori estimates} 
\label{sub:the_a_priori_estimates}


Define energy functionals and energy-dissipation functionals, for pure spatial derivatives, that,
  \begin{align}\label{func:higher-deri-x}
    E_{s}(t) & = |u|_{\hx{s}}^2 + \nm{(f_A,\,f_B)}_{\hxq{s}}^2 = E_{s,u} + E_{s,f}, \\[5pt]
    D_{s}(t) & = |\nabla_x u|_{\hx{s}}^2 + \sum_{k=0}^s \left[ \iint \abs{ \nabla_q \left( \tfrac{\nabla_x^k f_A}{\sqrt{M_A}} \right) }^2 M_A \d q \d x
            + \iint \abs{ \nabla_q \left( \tfrac{\nabla_x^k f_B}{\sqrt{M_B}} \right) }^2 M_B \d q \d x \right]\\\no
        & \quad + \nm{f_A - 2 f_B \sqrt{M_B}}^2_{\hxq{s}}
        \\\no & = D_{s,u} + D_{s,f} + D_{s,f,r},
  \end{align}
for spatial derivatives on first-order moment $j = \agl{qf}$ with index $s' = s-1$, that,
  \begin{align}\label{func:higher-deri-j}
    E_{s',j}(t) & = \nm{q (f_A,\,f_B)}_{\hxq{{s'}}}^2, \\[5pt]
    D_{s',j}(t) & = \sum_{k=0}^{s'} \left[
                        \iint \abs{ q \nabla_q \left( \tfrac{\nabla_x^k f_A}{\sqrt{M_A}} \right) }^2 M_A \d q \d x
                        + \iint \abs{ q \nabla_q \left( \tfrac{\nabla_x^k f_B}{\sqrt{M_B}} \right) }^2 M_B \d q \d x \right] \no\\
                & \quad + \nm{f_A - 2 f_B \sqrt{M_B}}^2_{\hxq{s'}},
  \end{align}
and, for mixed spatial-configurational derivatives, that,
  \begin{align}\label{func:higher-deri-mix}
    E_{s,\mix}(t) & = \sum_{\substack{k+l=s\\l\ge 1}} \eta^l \nm{\nabla_q^l (f_A,\,f_B)}_{\hxq{k}}^2 , \\[7pt]
    D_{s,\mix}(t) & = \sum_{\substack{k+l\le s\\[1pt] l\ge 1}} \eta^l \left[
            \iint \abs{ \nabla_q \big( \tfrac{\nabla_l^k f_A}{\sqrt{M_A}} \big) }^2 M_A \d q \d x
            + \iint \abs{ \nabla_q \big( \tfrac{\nabla_l^k f_B}{\sqrt{M_B}} \big) }^2 M_B \d q \d x
          \right] \no\\
        & \quad + \sum_{\substack{k+l\le s\\[1pt] l\ge 1}} \eta^l \nm{\nabla_l^k f_A - 2 \nabla_l^k f_B \sqrt{M_B}}_{\lxq}^2,
  \end{align}

By introducing the similar notations for number densities that
  \begin{align}
    E_{s,\rho}(t) = |\rho_A|_{\hx{s}}^2, \quad D_{s,\rho}(t) = |\rho_A|_{\hx{s}}^2,
  \end{align}
we will prove the following \emph{a-priori} estimate, as follows:

\begin{proposition}[A-priori estimate] \label{prop:a-priori}
  Denote that
  \begin{align}
    \E_s^\eta & = E_{s} + \eta E_{s,\rho} + \eta E_{s,\mix} + \eta^2 E_{s',j}, \\
    \D_s^\eta & = D_{s} + \eta D_{s,\rho} + \eta D_{s,\mix} + \eta^2 D_{s',j}.
  \end{align}

  Assuming that $\E_s^\eta \le \eps$, we can get, for some small fixed constant $\eta$, that
  \begin{align}\label{esm:a-priori}
    \frac{\d}{\d t} \E_s^\eta + \D_s^\eta
    \ls ( \eta^\frac{1}{2} + \eta^{-\frac{s}{2}} \eps^\frac{1}{2} ) \D_s^\eta .
  \end{align}
\end{proposition}

Note that we have the equivalence between the energy functionals defined in a way with or without the small constant parameter $\eta$, namely, $\E_s^\eta \sim \E_s$. Indeed, it is an easy matter to check $\eta^{-(s+1)} \E_s \le \E_s^\eta \le 2 \E_s$. Similar equivalence holds for dissipation functionals $\D_s^\eta \sim \D_s$.

\subsection{Higher-order spatial estimates for coupled system} 
\label{sub:spatial_estimates}


We will perform energy estimates of higher-order derivatives over variables $x$ on fluid velocity equation \eqref{esm:pure-x-u}, and microscopic equations of two species $(A,\,B)$ \eqref{esm:pure-x-AB}, respectively, then combine them together to give the pure spatial estimates \eqref{esm:pure-x-ABu}.

\smallskip\noindent\underline{\textbf{$H^s_x$-estimate for fluid velocity $u$}}:
Applying higher-order derivative operator $\nabla_x^s$ with index $s\ge 5$ on the third equation of \eqref{sys:perturbative}, we have
  \begin{multline}\label{esm:pure-x-u}
    \frac{1}{2} \frac{\d}{\d t} |\nabla_x^s u|_{\lx}^2 + \mu |\nabla_x^{s+1} u|_{\lx}^2
    \ls \abs{u}_{\hx{s}} \abs{\nabla_x u}_{\hx{s-1}} \abs{\nabla_x^s u}_{\lx}  \\
      + \skpa{\nabla_q U_A \otimes q \nabla_x^s f_A \sqrt{M_A}}{\nabla_x^{s+1} u}
      + \skpa{\nabla_q U_B \otimes q \nabla_x^s f_B \sqrt{M_B}}{\nabla_x^{s+1} u}.
  \end{multline}
The process is direct, by noticing the
  \begin{align*}
    \skpa{\nabla_x^s (u \cdot \nabla_x u)}{\nabla_x^s u}
    & = \skpa{[\nabla_x^s,\, u \cdot \nabla_x]u}{\nabla_x^s u} + \skpa{u \cdot \nabla_x (\nabla_x^s u)}{\nabla_x^s u} \\\no
    & \ls \abs{u}_{\hx{s}} \abs{\nabla_x u}_{\hx{s-1}} \cdot \abs{\nabla_x^s u}_{\lx}.
  \end{align*}
where the latter term in first line vanishes due to the incompressibility condition $\div_x u =0$, while the commutator term is bounded by using the Moser-type inequality (see \cite{Maj-84b,Tay-PDE3}):
  \begin{align*}
    \abs{[\nabla_x^s,\, u \cdot \nabla_x]v}_{\lx}
    \ls \abs{u}_{\hx{s}} \abs{\nabla_x v}_{L^\infty_x}
        + \abs{\nabla_x u}_{L^\infty_x} \abs{\nabla_x v}_{\hx{s-1}},
  \end{align*}
and the Sobolev embedding inequalities together.

\smallskip\noindent\underline{\textbf{$H^s_x L^2_q$-estimate for microscopic perturbations $(f_A,\, f_B)$}}: Apply higher-order spatial derivative operator $\nabla_x^s$ with index $s\ge 5$ on the first and second equation of \eqref{sys:perturbative}, we then get, respectively,
  \begin{multline}\label{Eq:deri-x-A}
    \p_t \nabla_x^s f_A + \nabla_x^s (u \cdot \nabla_x f_A) + \nabla_x^s (\nabla_x u q \nabla_q f_A) = \mathcal{L}_A (\nabla_x^s f_A) \\ + \nabla_x^s [\nabla_x u q \nabla_q U_A(\sqrt{M_A} + \tfrac{1}{2} f_A)] + \nabla_x^s r_A ,
  \end{multline}
and
  \begin{multline}\label{Eq:deri-x-B}
    \p_t \nabla_x^s f_B + \nabla_x^s (u \cdot \nabla_x f_B) + \nabla_x^s (\nabla_x u q \nabla_q f_B) = \mathcal{L}_B (\nabla_x^s f_B) \\ + \nabla_x^s [\nabla_x u q \nabla_q U_B (\sqrt{M_B} + \tfrac{1}{2} f_B)] + \nabla_x^s r_B .
  \end{multline}
Here we actually have used the commutative property between the spatial derivative operator $\nabla_x^s$ and the Fokker-Planck operator $\L_{\alpha}$ with $\alpha=A,\, B$.

Take an inner product with $(\nabla_x^s f_A,\, \nabla_x^s f_B)$, respectively, in space $L^2_{x,q}$, then we can get some similar estimates for above two equations, except the last reaction terms. We take species $A$ for example, and write the estimates as follows:
  \begin{align}
    \skpa{\p_t \nabla_x^s f_A}{\nabla_x^s f_A} = \frac{1}{2} \frac{\d}{\d t} \nm{\nabla_x^s f_A}_{\lxq}^2,
  \end{align}
and
  \begin{align}
    \skpa{\nabla_x^s (u \cdot \nabla_x f_A)}{\nabla_x^s f_A}
    = \skpa{[\nabla_x^s,\, u \cdot \nabla_x]f_A}{\nabla_x^s f_A} + \skpa{u \cdot \nabla_x (\nabla_x^s f_A)}{\nabla_x^s f_A},
  \end{align}
where the latter product will vanish due to the incompressibility condition $\div_x u =0$. The commutator in the former can be decomposed as:
  \begin{align}
    \nm{[\nabla_x^s,\, u \cdot \nabla_x]f_A}_{\lxq}
    = \sum_{\substack{s_1 + s_2 = s\\s_1 \ge 1}} \nm{\nabla_x^{s_1} u \nabla_x (\nabla_x^{s_2} f_A)}_{\lxq}.
  \end{align}

In the case $s_1 = s$, the H\"older inequality yields that,
  \begin{align}
    \nm{\nabla_x^{s} u \nabla_x f_A}_{\lxq} \le \abs{\nabla_x^s u}_{\lx} \nm{\nabla_x f_A}_{L^\infty_x L^2_q}
    \ls \abs{\nabla_x^s u}_{\lx} \nm{\nabla_x f_A}_{H^2_x L^2_q},
  \end{align}
where we have used the Sobolev embedding inequality in the last line.

In the case $s_1 = s-1$, we get similarly, that
  \begin{align}
    \nm{\nabla_x^{s-1} u (\nabla_x^{2} f_A)}_{\lxq}
    \le \abs{\nabla_x^{s-1} u}_{L^4_x} \nm{\nabla_x^2 f_A}_{L^4_x L^2_q}
    \ls \abs{\nabla_x^{s-1} u}_{H^1_x} \nm{\nabla_x^2 f_A}_{H^1_x L^2_q},
  \end{align}
where we have used the H\"older inequality and the Sobolev embedding inequality again.

For the left cases $1 \le s_1 \le s-2$ (and hence $s_2 \le s-1$), we can get the following control,
  \begin{align}
    \sum_{\substack{s_1 + s_2 = s\\1\le s_1 \le s-2}} \nm{\nabla_x^{s_1} u \nabla_x (\nabla_x^{s_2} f_A)}_{\lxq}
    \le \sum_{\substack{s_1 + s_2 = s\\1\le s_1 \le s-2}} \abs{\nabla_x^{s_1} u}_{L^\infty_x} \nm{\nabla_x^{s_2+1} f_A}_{L^2_x L^2_q}
    \ls \abs{u}_{H^s_x} \nm{f_A}_{H^s_x L^2_q}.
  \end{align}

Combining three above cases together gives the bound on the commutator
  \begin{align}
    \nm{[\nabla_x^s,\, u \cdot \nabla_x]f_A}_{\lxq} \ls \abs{u}_{H^s_x} \nm{f_A}_{H^s_x L^2_q},
  \end{align}
and consequently gives the control on the term:
  \begin{align}
    \skpa{\nabla_x^s (u \cdot \nabla_x f_A)}{\nabla_x^s f_A}
    \ls \abs{u}_{H^s_x} \nm{f_A}_{H^s_x L^2_q} \nm{\nabla_x^{s} f_A}_{\lxq}.
  \end{align}

For the third term in left-hand side of equation \eqref{Eq:deri-x-A},  by noticing the fact $\nabla_q \cdot (\nabla_x u q) = \div_x u =0$, we have the following decomposition,
  \begin{align}
    \skpa{\nabla_x^s (\nabla_x u q \nabla_q f_A)}{\nabla_x^s f_A}
    & = \skpa{\nabla_x^s (\nabla_x u q \nabla_q f_A) - \nabla_x u q \nabla_q (\nabla_x^s f_A)}{\nabla_x^s f_A} \\\no
    & = \sum_{\substack{s_1+s_2=s\\s_1 \ge 1}} \skpa{\nabla_x^{s_1} (\nabla_x u) \nabla_q \nabla_x^{s_2} f_A}{q \nabla_x^s f_A}.
  \end{align}
The summation can be bounded in a similar discussion way by cases as above, that,
  \begin{align}
    & \sum_{\substack{s_1+s_2=s\\s_1 \ge 1}} \nm{\nabla_x^{s_1} (\nabla_x u) \nabla_q \nabla_x^{s_2} f_A}_{\lxq} \\\no
    & = \nm{\nabla_x^{s+1} u \nabla_q f_A}_{\lxq}
        + \nm{\nabla_x^{s} u \nabla_q \nabla_x f_A}_{\lxq}
        + \sum_{\substack{s_1+s_2=s\\1\le s_1 \le s-2}} \nm{\nabla_x^{s_1+1} u \nabla_q \nabla_x^{s_2} f_A}_{\lxq} \\\no
    & \le \abs{\nabla_x^{s+1} u}_{\lx} \nm{\nabla_q f_A}_{L^\infty_x L^2_q}
          + \abs{\nabla_x^{s} u}_{L^4_x} \nm{\nabla_q \nabla_x f_A}_{L^4_x L^2_q}
          + \sum_{\substack{s_1+s_2=s\\1\le s_1 \le s-2}} \abs{\nabla_x^{s_1+1} u}_{L^\infty_x} \nm{\nabla_q \nabla_x^{s_2} f_A}_{L^2_x L^2_q} \\\no
    & \ls \abs{\nabla_x^{s+1} u}_{\lx} \nm{\nabla_q f_A}_{H^2_x L^2_q}
          + \abs{\nabla_x^{s} u}_{H^1_x} \nm{\nabla_q \nabla_x f_A}_{H^1_x L^2_q}
          + \abs{\nabla_x u}_{\hx{s}} \nm{\nabla_q f_A}_{\hxq{s-1}} \\\no
    & \ls \abs{\nabla_x u}_{\hx{s}} \nm{\nabla_q f_A}_{\hxq{s-1}},
  \end{align}
where we have used again the H\"older inequality and the Sobolev embedding inequalities. Hence, we get,
  \begin{align}
    \skpa{\nabla_x^s (\nabla_x u q \nabla_q f_A)}{\nabla_x^s f_A}
    \ls \abs{\nabla_x u}_{\hx{s}} \nm{\nabla_q f_A}_{\hxq{s-1}} \nm{\nabla_q U_A \nabla_x^s f_A}_{\lxq}.
  \end{align}

We now turn to consider the right-hand side of equation \eqref{Eq:deri-x-A}. Firstly we have,
  \begin{align}
    \skpa{\L_A \nabla_x^s f_A}{\nabla_x^s f_A}
    = - \iint \abs{ \nabla_q \left( \tfrac{\nabla_x^s f_A}{\sqrt{M_A}} \right) }^2 M_A \d q \d x
    = - \nm{\nabla_q ( \tfrac{\nabla_x^s f_A}{\sqrt{M_A}} )}_M^2.
  \end{align}

For the penultimate term in the right-hand side, we get
  \begin{align}
    & \skpa{\nabla_x^s (\nabla_x u q \nabla_q U_A f_A)}{\nabla_x^s f_A} \\\no
    & = \sum_{s_1+s_2=s} \skpa{\nabla_x^{s_1+1} u q \nabla_x^{s_2} f_A}{\nabla_q U_A \nabla_x^s f_A} \\\no
    & \ls (\abs{\nabla_x u}_{\hx{s}} \nm{q f_A}_{\hxq{2}} + \abs{u}_{\hx{s}} \nm{q f_A}_{\hxq{s}}) \nm{\nabla_q U_A \nabla_x^s f_A}_{\lxq},
  \end{align}
where we have used a similar discussion by cases as above, that
  \begin{align}
    & \sum_{s_1+s_2=s} \nm{\nabla_x^{s_1+1} u q \nabla_x^{s_2} f_A}_{\lxq} \\\no
    & = \nm{\nabla_x^{s+1} u q f_A}_{\lxq} + \nm{\nabla_x^{s} u q \nabla_x f_A}_{\lxq}
      + \nm{\nabla_x^{s-1} u q \nabla_x^{2} f_A}_{\lxq}
      + \sum_{s_1 \le s-3} \nm{\nabla_x^{s_1+1} u q \nabla_x^{s_2} f_A}_{\lxq} \\\no
    & \le \abs{\nabla_x^{s+1} u}_{\lx} \nm{q f_A}_{L^\infty_x L^2_q}
          + \abs{\nabla_x^{s} u}_{L^2_x} \nm{q \nabla_x f_A}_{L^\infty_x L^2_q}
          + \abs{\nabla_x^{s-1} u}_{L^4_x} \nm{q \nabla_x^2 f_A}_{L^4_x L^2_q} \\\no
        & \quad
          + \sum_{s_1 \le s-3} \abs{\nabla_x^{s_1+1} u}_{L^\infty_x} \nm{q \nabla_x^{s_2} f_A}_{L^2_x L^2_q} \\\no
    & \ls \abs{\nabla_x u}_{\hx{s}} \nm{q f_A}_{\hxq{2}} + \abs{u}_{\hx{s}} \nm{q f_A}_{\hxq{s}}.
  \end{align}
Note that the H\"older inequality and the Sobolev embedding inequalities are used again.

Since the estimates for equation \eqref{Eq:deri-x-B} of species $B$ are similar as above process for species $A$, then we are left to deal with the reaction contribution now, by combining both equations for two species $A,\, B$ together. We calculate that
  \begin{align}
    & \skpa{\nabla_x^s r_A}{\nabla_x^s f_A} + \skpa{\nabla_x^s r_B}{\nabla_x^s f_B} \\\no
    & = - \skpa{\nabla_x^s(f_A - 2f_B \sqrt{M_B} - f_B^2)}{\nabla_x^s f_A}
        + 2 \skpa{\sqrt{M_B} \nabla_x^s(f_A - 2f_B \sqrt{M_B} - f_B^2)}{\nabla_x^s f_B} \\\no
    & = - \nm{\nabla_x^s (f_A - 2f_B \sqrt{M_B})}_{\lxq}^2 + \iint \nabla_x^s (f_B^2) \cdot \nabla_x^s (f_A - 2f_B \sqrt{M_B}) \d q \d x \\\no
    & \gs - \nm{\nabla_x^s (f_A - 2f_B \sqrt{M_B})}_{\lxq}^2 + \nm{f_B}_{H^2_x H^2_q} \nm{f_B}_{\hxq{s}} \nm{\nabla_x^s (f_A - 2f_B \sqrt{M_B})}_{\lxq},
  \end{align}
where we have used the product formula in Sobolev spaces
  \begin{align*}
    \abs{\nabla_x^s (f g)}_{\lx} \ls \abs{f}_{L^\infty_x} \abs{g}_{\hx{s}} + \abs{f}_{\hx{s}} \abs{g}_{L^\infty_x}.
  \end{align*}

Therefore, combining all the above estimates together enable us to get the higher-order purely spatial estimates on microscopic perturbations $(f_A,\, f_B)$, as follows,
  \begin{align}\label{esm:pure-x-AB}
    & \frac{1}{2} \frac{\d}{\d t} \nm{\nabla_x^s (f_A,\, f_B)}_{\lxq}^2
      + \nm{\nabla_q ( \tfrac{\nabla_x^s f_A}{\sqrt{M_A}} )}_M^2
              + \nm{\nabla_q ( \tfrac{\nabla_x^s f_B}{\sqrt{M_B}} )}_M^2
              + \nm{\nabla_x^s (f_A - 2 f_B \sqrt{M_B})}_{\lxq}^2
  \no\\
     & \ls \skpa{\nabla_x^{s+1} u q \nabla_q U_A \sqrt{M_A}}{\nabla_x^s f_A} + \skpa{\nabla_x^{s+1} u q \nabla_q U_B \sqrt{M_B}}{\nabla_x^s f_B}
     \\[5pt] \no
        & \quad + |u|_{\hx{s}} \nm{f_A}_{\hxq{s}} \nm{\nabla_x^s f_A}_{\lxq}
          + |\nabla_x u|_{\hx{s}} \nm{\nabla_q f_A}_{\hxq{s-1}} \nm{\nabla_q U_A \nabla_x^s f_A}_{\lxq} \\\no
        & \quad + \left( |u|_{\hx{s}} \nm{qf_A}_{\hxq{s}} + |\nabla_x u|_{\hx{s}} \nm{q f_A}_{\hxq{2}} \right) \nm{\nabla_q U_A \nabla_x^s f_A}_{\lxq} \\[4pt] \no
        & \quad + |u|_{\hx{s}} \nm{f_B}_{\hxq{s}} \nm{\nabla_x^s f_B}_{\lxq}
                + |\nabla_x u|_{\hx{s}} \nm{\nabla_q f_B}_{\hxq{s-1}} \nm{\nabla_q U_B \nabla_x^s f_B}_{\lxq} \\[3pt] \no
        & \quad + \left( |u|_{\hx{s}} \nm{qf_B}_{\hxq{s}} + |\nabla_x u|_{\hx{s}} \nm{q f_B}_{\hxq{2}} \right) \nm{\nabla_q U_B \nabla_x^s f_B}_{\lxq} \\[4pt] \no
        & \quad + \nm{f_B}_{H^2_x H^2_q} \nm{f_B}_{\hxq{s}} \nm{\nabla_x^s (f_A - 2 f_B \sqrt{M_B})}_{\lxq}.
  \end{align}

\smallskip\noindent\underline{\textbf{Higher-order spatial estimates for coupled system}}: At the end of this subsection, by noticing the cancellation relations between micro-macro coupling:
  \begin{align}\label{eq:cancel-x}
  \begin{cases}
    \skpa{\nabla_q U_A \otimes q \nabla_x^s f_A \sqrt{M_A}}{\nabla_x^{s+1} u}
      + \skpa{\nabla_x^{s+1} u q \nabla_q U_A \sqrt{M_A}}{\nabla_x^s f_A} = 0,
    \\[4pt]
    \skpa{\nabla_q U_B \otimes q \nabla_x^s f_B \sqrt{M_B}}{\nabla_x^{s+1} u}
      + \skpa{\nabla_x^{s+1} u q \nabla_q U_B \sqrt{M_B}}{\nabla_x^s f_B} = 0,
  \end{cases}
  \end{align}
we can conclude the higher-order spatial estimates for coupled perturbative system \eqref{sys:perturbative}, that
  \begin{align}\label{esm:pure-x-ABu}
     & \frac{1}{2} \frac{\d}{\d t} \left( |\nabla_x^s u|_{\lx}^2 + \nm{\nabla_x^s (f_A,\,f_B)}_{\lxq}^2 \right)
        + \mu |\nabla_x^{s+1} u|_{\lx}^2 \\\no & \quad
        + \nm{\nabla_q ( \tfrac{\nabla_x^s f_A}{\sqrt{M_A}} )}_M^2
        + \nm{\nabla_q ( \tfrac{\nabla_x^s f_B}{\sqrt{M_B}} )}_M^2
        + \nm{\nabla_x^s (f_A - 2 f_B \sqrt{M_B})}_{\lxq}^2
  \\\no
     & \ls |u|_{\hx{s}} |\nabla_x u|_{\hx{s}} |\nabla_x^s u|_{\lx}
            + |u|_{\hx{s}} \nm{f_A}_{\hxq{s}} \nm{\nabla_x^s f_A}_{\lxq}
            + |\nabla_x u|_{\hx{s}} \nm{\nabla_q f_A}_{\hxq{s-1}} \nm{\nabla_q U_A \nabla_x^s f_A}_{\lxq} \\\no
        & \quad + \left( |u|_{\hx{s}} \nm{qf_A}_{\hxq{s}} + |\nabla_x u|_{\hx{s}} \nm{q f_A}_{\hxq{2}} \right) \nm{\nabla_q U_A \nabla_x^s f_A}_{\lxq} \\[4pt] \no
        & \quad + |u|_{\hx{s}} \nm{f_B}_{\hxq{s}} \nm{\nabla_x^s f_B}_{\lxq}
               + |\nabla_x u|_{\hx{s}} \nm{\nabla_q f_B}_{\hxq{s-1}} \nm{\nabla_q U_B \nabla_x^s f_B}_{\lxq} \\[3pt] \no
        & \quad  + \left( |u|_{\hx{s}} \nm{qf_B}_{\hxq{s}} + |\nabla_x u|_{\hx{s}} \nm{q f_B}_{\hxq{2}} \right) \nm{\nabla_q U_B \nabla_x^s f_B}_{\lxq} \\[4pt] \no
        & \quad + \nm{f_B}_{H^2_x H^2_q} \nm{f_B}_{\hxq{s}} \nm{\nabla_x^s (f_A - 2 f_B \sqrt{M_B})}_{\lxq}.
  \end{align}

By virtue of notations for energy and dissipation functionals defined in subsection \S \ref{sub:the_a_priori_estimates}, we can write
  \begin{align}
    \frac{\d}{\d t} E_{s} + D_{s} 
     & \ls E_{s,u}^\frac{1}{2} D_{s,u} + E_{s,u}^\frac{1}{2} (D_{s,f} + D_{s,\rho}) 
          + D_{s,u}^\frac{1}{2} \nm{\nabla_q (f_A,\,f_B)}_{\hxq{s-1}} (D_{s,f}^\frac{1}{2} + D_{s,\rho}^\frac{1}{2}) \\[3pt] \no
     & \quad + \left[ E_{s,u}^\frac{1}{2} \nm{q(f_A,\,f_B)}_{\hxq{s}} + D_{s,u}^\frac{1}{2} \nm{q(f_A,\,f_B)}_{\hxq{2}} \right]  (D_{s,f}^\frac{1}{2} + D_{s,\rho}^\frac{1}{2}) \\[3pt] \no
     & \quad + \nm{f_B}_{H^2_x H^2_q} \nm{f_B}_{\hxq{s}} D_{s,f,r}^\frac{1}{2}.
  \end{align}

\begin{remark}\label{rem:esm-need}
  Indeed, this suggests us to estimate higher-order mixed derivatives for $(f_A,\,f_B)$ and higher-order purely spatial derivatives for first-order moments $(qf_A,\,qf_B)$. In fact, for the latter, it suffices to control higher but not highest order derivatives, since the highest order quantities should be viewed as dissipative part.
\end{remark}

\subsection{Higher-order spatial estimates on first-order moments} 
\label{sub:first_order_moments}


We will perform in this subsection higher-order spatial derivative estimates for first-order moments $(qf_A,\,qf_B)$ with a lower Sobolev index $s'=s-1$ (actually, $4\le s' \le s-1$ suffices to get our result).

We consider to act higher-order spatial derivative operator $\nabla_x^k$ with index $4 \le k \le s'$ on the first and second equation of \eqref{sys:perturbative}, then by taking inner product of equation \eqref{Eq:deri-x-A} with $|q|^2 \nabla_x^k f_A$, we get,
  \begin{align}
    \skpa{\p_t \nabla_x^k f_A}{|q|^2 \nabla_x^k f_A} = \frac{1}{2} \frac{\d}{\d t} \nm{q \nabla_x^k f_A}_{\lxq}^2,
  \end{align}
and
  \begin{align}
    \skpa{\nabla_x^k (u \cdot \nabla_x f_A)}{|q|^2 \nabla_x^k f_A}
    & = \skpa{[\nabla_x^k,\, u \cdot \nabla_x]q f_A}{q \nabla_x^k f_A} + \skpa{u \cdot \nabla_x (q\nabla_x^k f_A)}{q \nabla_x^k f_A}, \\\no
    & \ls \abs{u}_{H^k_x} \nm{qf_A}_{H^k_x L^2_q} \nm{q\nabla_x^{k} f_A}_{\lxq},
  \end{align}
where we have used the same discussion as that in subsection \S\ref{sub:first_order_moments} and the incompressibility $\div_x u=0$ again.

We decompose the estimate for the third term in left-hand side into three parts, more precisely,
  \begin{align}
    & \skpa{\nabla_x^k (\nabla_x u q \nabla_q f_A)}{|q|^2 \nabla_x^k f_A} \\\no
    & = \skpa{\nabla_x^{k} (\nabla_x u) q \nabla_q f_A}{|q|^2 \nabla_x^k f_A}
        + \skpa{\nabla_x u q \nabla_q \nabla_x^{k} f_A}{|q|^2 \nabla_x^k f_A} \\\no
        & \quad + \sum_{\substack{k_1+k_2=k\\1 \le k_1 \le k-1}} \skpa{\nabla_x^{k_1} (\nabla_x u) q \nabla_q \nabla_x^{k_2} f_A}{|q|^2 \nabla_x^k f_A} .
  \end{align}
The first term in above expression can be controlled as follows,
  \begin{align}
    \skpa{\nabla_x^{k} (\nabla_x u) q \nabla_q f_A}{|q|^2 \nabla_x^k f_A}
      \ls \abs{\nabla_x u}_{\hx{k}} \nm{q \nabla_q f_A}_{\hxq{2}} \nm{q \nabla_q U_A \nabla_x^k f_A}_{\lxq}.
  \end{align}
Recalling the fact $\nabla_q \cdot (\nabla_x u q) = \div_x u =0$, we have,
  \begin{align}
    \skpa{\nabla_x u q \nabla_q \nabla_x^{k} f_A}{|q|^2 \nabla_x^k f_A}
    & = \iint \nabla_x u q |q|^2 \nabla_q [\tfrac{1}{2} (\nabla_x^{k} f_A)^2] \d q \d x \\\no
    & = - \iint \nabla_x u q \otimes q |\nabla_x^{k} f_A|^2 \d q \d x \\\no
    & \le \abs{\nabla_x u}_{L^\infty_x} \nm{q \nabla_x^k f_A}_{\lxq}^2.
  \end{align}
For the summation term, we get
  \begin{align}
    & \sum_{\substack{k_1+k_2=k\\1 \le k_1 \le k-1}} \skpa{\nabla_x^{k_1} (\nabla_x u) q \nabla_q \nabla_x^{k_2} f_A}{|q|^2 \nabla_x^k f_A} \\\no
    & = \skpa{\nabla_x^k u q \nabla_q \nabla_x f_A}{|q|^2 \nabla_x^k f_A}
        + \skpa{\nabla_x^{k-1} u q \nabla_q \nabla_x^{2} f_A}{|q|^2 \nabla_x^k f_A} \\\no
        & \quad + \sum_{\substack{k_1+k_2=k\\1 \le k_1 \le k-3}} \skpa{\nabla_x^{k_1} (\nabla_x u) q \nabla_q \nabla_x^{k_2} f_A}{|q|^2 \nabla_x^k f_A} \nm{|q|^2 \nabla_x^k f_A}_{\lxq}
    \\\no
    & \le \abs{\nabla_x^k u}_{L^2_x} \nm{q \nabla_q \nabla_x f_A}_{L^\infty_x L^2_q} \nm{|q|^2 \nabla_x^k f_A}_{\lxq}
        + \abs{\nabla_x^{k-1} u}_{L^4_x} \nm{q \nabla_q \nabla_x^2 f_A}_{L^4_x L^2_q} \nm{|q|^2 \nabla_x^k f_A}_{\lxq} \\\no
        & \quad + \sum_{\substack{k_1+k_2=k\\1 \le k_1 \le k-3}} \abs{\nabla_x^{k_1 +1} u}_{L^\infty_x} \nm{q \nabla_q \nabla_x^{k_2} f_A}_{L^2_x L^2_q} \nm{|q|^2 \nabla_x^k f_A}_{\lxq}
    \\\no
    & \ls \abs{u}_{\hx{k}} \nm{q \nabla_q f_A}_{\hxq{k-1}} \nm{q \nabla_q U_A \nabla_x^k f_A}_{\lxq}.
  \end{align}

As a consequence, we obtain that,
  \begin{multline}
    \skpa{\nabla_x^k (\nabla_x u q \nabla_q f_A)}{|q|^2 \nabla_x^k f_A}
    \ls \abs{u}_{\hx{3}} \nm{q \nabla_x^k f_A}_{\lxq}^2 \\
        + (\abs{\nabla_x u}_{\hx{k}} \nm{q \nabla_q f_A}_{\hxq{2}} + \abs{u}_{\hx{k}} \nm{q \nabla_q f_A}_{\hxq{k-1}}) \nm{q \nabla_q U_A \nabla_x^k f_A}_{\lxq}.
  \end{multline}

Concerning the terms on the right-hand side, we derive from the expression of Fokker-Planck operator that,
  \begin{align}
    \skpa{\L_A \nabla_x^k f_A}{|q|^2 \nabla_x^k f_A}
    & = -\skpa{\nabla_q (\tfrac{\nabla_x^k f_A}{\sqrt{M_A}}) }{\nabla_q (|q|^2 \tfrac{\nabla_x^k f_A}{\sqrt{M_A}})} \\\no
    & = - \iint \abs{ q \nabla_q (\tfrac{\nabla_x^k f_A}{\sqrt{M_A}}) }^2 M_A \d q \d x
        - \iint q M_A \nabla_q \big( \abs{\tfrac{\nabla_x^k f_A}{\sqrt{M_A}}}^2 \big) \d q \d x \\\no
    & = - \iint \abs{ q \nabla_q (\tfrac{\nabla_x^k f_A}{\sqrt{M_A}}) }^2 M_A \d q \d x
        + 3 \iint \abs{\tfrac{\nabla_x^k f_A}{\sqrt{M_A}}}^2 M_A \d q \d x \\\no & \quad
        - \iint q \nabla_q U_A \abs{\tfrac{\nabla_x^k f_A}{\sqrt{M_A}}}^2 M_A \d q \d x \\\no
    & = - \nm{q \nabla_q ( \tfrac{\nabla_x^k f_A}{\sqrt{M_A}} )}_M^2 + 3 \nm{\nabla_x^k f_A}_{\lxq}^2
        - \iint q \nabla_q U_A \abs{\nabla_x^k f_A}^2 \d q \d x.
  \end{align}
Note here the third term in the last line can be controlled by the H\"older inequality and Lemma \ref{lemm:Poinc-inequ-0} as follows,
  \begin{align}
    \iint q \nabla_q U_A \abs{\nabla_x^k f_A}^2 \d q \d x
    & \le \nm{q \nabla_q U_A \nabla_x^k f_A}_{\lxq} \nm{\nabla_x^k f_A}_{\lxq} \\\no
    & \ls \big(\nm{\agl{q} \nabla_q ( \tfrac{\nabla_x^k f_A}{\sqrt{M_A}})}_M + \abs{\nabla_x^k \rho_A}_{\lx} \big) \cdot \big(\nm{\nabla_q ( \tfrac{\nabla_x^k f_A}{\sqrt{M_A}})}_M + \abs{\nabla_x^k \rho_A}_{\lx} \big) \\\no
    & \le \tfrac{1}{4} \nm{q \nabla_q ( \tfrac{\nabla_x^k f_A}{\sqrt{M_A}})}_M^2 + C \big( \nm{\nabla_q ( \tfrac{\nabla_x^k f_A}{\sqrt{M_A}})}_M + \abs{\nabla_x^k \rho_A}_{\lx} \big)^2,
  \end{align}
where we have used the simple fact $\abs{\nabla_x^k \rho_A}_{\lx} \ls \nm{\nabla_x^k f_A}_{\lxq}$ due to H\"older's inequality again.
Thus, we are ready to obtain the estimates,
  \begin{align}
    \skpa{\L_A \nabla_x^k f_A}{|q|^2 \nabla_x^k f_A}
    \ge - \tfrac{3}{4} \nm{q \nabla_q ( \tfrac{\nabla_x^k f_A}{\sqrt{M_A}} )}_M^2
        + C \big( \nm{\nabla_q ( \tfrac{\nabla_x^k f_A}{\sqrt{M_A}})}_M^2 + \abs{\nabla_x^k \rho_A}_{\lx}^2 \big).
  \end{align}

As for the penultimate term in the right-hand side, we can infer from a discussion by cases that,
  \begin{align}
    & \skpa{\nabla_x^k (\nabla_x u q \nabla_q U_A f_A)}{|q|^2 \nabla_x^k f_A} \\\no
    & = \sum_{k_1+k_2=k} \skpa{\nabla_x^{k_1+1} u q \nabla_q U_A \nabla_x^{k_2} f_A}{|q|^2 \nabla_x^k f_A} \\\no
    & \ls (\abs{\nabla_x u}_{\hx{k}} \nm{q \nabla_q U_A f_A}_{\hxq{2}} + \abs{u}_{\hx{k}} \nm{q \nabla_q U_A f_A}_{\hxq{k}}) \nm{q \nabla_q U_A \nabla_x^k f_A}_{\lxq}.
  \end{align}
and
  \begin{align}
    \skpa{\nabla_x^k (\nabla_x u q \nabla_q U_A \sqrt{M_A})}{|q|^2 \nabla_x^k f_A}
    \le \abs{\nabla_x u}_{\hx{k}} \nm{|q|^2 \sqrt{M_A}}_{L^2_q} \nm{q \nabla_q U_A \nabla_x^k f_A}_{\lxq}. 
  \end{align}

We next deal with the reaction contribution. Combining the two equations for both species $A,\, B$ together implies,
  \begin{align}
    & \skpa{\nabla_x^k r_A}{|q|^2 \nabla_x^k f_A} + \skpa{\nabla_x^k r_B}{|q|^2 \nabla_x^k f_B} \\\no
    & = - \nm{q \nabla_x^k (f_A - 2f_B \sqrt{M_B})}_{\lxq}^2 + \iint \nabla_x^k (|q|^2 f_B^2) \cdot \nabla_x^k (f_A - 2f_B \sqrt{M_B}) \d q \d x \\\no
    & \ls - \nm{q \nabla_x^k (f_A - 2f_B \sqrt{M_B})}_{\lxq}^2 + \nm{qf_B}_{H^2_x H^2_q} \nm{qf_B}_{\hxq{k}} \nm{\nabla_x^k (f_A - 2f_B \sqrt{M_B})}_{\lxq}.
  \end{align}
Therefore, together with previous estimates for species $A$ and $B$, this concludes spatial estimates on first order moments $(qf_A,\,qf_B)$ with index $4 \le k \le s'=s-1$,
  \begin{align}\label{esm:pure-x-qAqB}
     & \frac{1}{2} \frac{\d}{\d t} \nm{q \nabla_x^k (f_A,\, f_B)}_{\lxq}^2
     + \tfrac{3}{4} \Big[ \nm{q \nabla_q ( \tfrac{\nabla_x^k f_A}{\sqrt{M_A}} )}_M^2 + \nm{q \nabla_q ( \tfrac{\nabla_x^k f_B}{\sqrt{M_B}} )}_M^2 \Big]
     + \nm{q\nabla_x^k (f_A - 2 f_B \sqrt{M_B})}_{\lxq}^2 \no\\
     & \ls \big( \nm{\nabla_q ( \tfrac{\nabla_x^k f_A}{\sqrt{M_A}})}_M^2 + \nm{\nabla_q ( \tfrac{\nabla_x^k f_B}{\sqrt{M_B}} )}_M^2 \big) + \abs{\nabla_x^k (\rho_A,\, \rho_B)}_{\lx}^2 \\\no
     & \quad + |u|_{\hx{k}} \nm{q f_A}_{\hxq{k}} \nm{q \nabla_x^k f_A}_{\lxq}
             + |u|_{\hx{3}} \nm{q \nabla_x^k f_A}_{\lxq}^2 \\\no
     & \quad + \left( |u|_{\hx{k}} \nm{q \nabla_q f_A}_{\hxq{k-1}} + |\nabla_x u|_{\hx{k}} \nm{q \nabla_q f_A}_{\hxq{2}} \right) \nm{q \nabla_q U_A \nabla_x^k f_A}_{\lxq} \\[4pt] \no
     & \quad + \left( |u|_{\hx{k}} \nm{q \nabla_q U_A f_A}_{\hxq{k}} + |\nabla_x u|_{\hx{k}} \nm{q \nabla_q U_A f_A}_{\hxq{2}} \right) \nm{q \nabla_q U_A \nabla_x^k f_A}_{\lxq} \\ \no
     & \quad + |\nabla_x u|_{\hx{k}} \nm{q \nabla_q U_A \nabla_x^k f_A}_{\lxq}
             + |\nabla_x u|_{\hx{k}} \nm{q \nabla_q U_B \nabla_x^k f_B}_{\lxq}
  \\\no
     & \quad + |u|_{\hx{k}} \nm{q f_B}_{\hxq{k}} \nm{q \nabla_x^k f_B}_{\lxq}
             + |u|_{\hx{3}} \nm{q \nabla_x^k f_B}_{\lxq}^2 \\\no
     & \quad + \left( |u|_{\hx{k}} \nm{q \nabla_q f_B}_{\hxq{k-1}} + |\nabla_x u|_{\hx{k}} \nm{q \nabla_q f_B}_{\hxq{2}} \right) \nm{q \nabla_q U_B \nabla_x^k f_B}_{\lxq} \\[4pt] \no
     & \quad + \left( |u|_{\hx{k}} \nm{q \nabla_q U_B f_B}_{\hxq{k}} + |\nabla_x u|_{\hx{k}} \nm{q \nabla_q U_B f_B}_{\hxq{2}} \right) \nm{q \nabla_q U_B \nabla_x^k f_B}_{\lxq}
  \\\no
     & \quad + \nm{q f_B}_{\hxq{k}} \nm{q f_B}_{H^2_x H^2_q} \nm{\nabla_x^k (f_A - 2 f_B \sqrt{M_B})}_{\lxq},
  \end{align}
whence we get
  \begin{align}
    \frac{\d}{\d t} E_{s',j} + & D_{s',j} 
      \ls (D_{s',f} + D_{s',\rho})
         + E_{s',u}^\frac{1}{2} (D_{s',f} + D_{s',\rho}) \\\no
      & + [
                   E_{s',u}^\frac{1}{2} \nm{q \nabla_q (f_A,\,f_B)}_{\hxq{s'-1}}
                   + E_{s,u}^\frac{1}{2} \nm{q \nabla_q (f_A,\,f_B)}_{\hxq{2}}
                ] (D_{s',j}^\frac{1}{2} + D_{s',f}^\frac{1}{2} + D_{s',\rho}^\frac{1}{2})
      \\[3pt] \no
      & + (E_{s',u}^\frac{1}{2} + E_{s,u}^\frac{1}{2}) (D_{s',j} + D_{s',f} + D_{s',\rho})
              + D_{s',u}^\frac{1}{2} (D_{s',j}^\frac{1}{2} + D_{s',f}^\frac{1}{2} + D_{s',\rho}^\frac{1}{2})
             \\\no
      & + \nm{q f_B}_{\hxq{s'}} \nm{q f_B}_{H^2_x H^2_q} D_{s',f,r}^\frac{1}{2},
  \end{align}
where we have used $D_{s',u} \le E_{s,u}$ and $\nm{q \nabla_q U_\alpha \nabla_x^k f_\alpha}_{\lxq} \ls D_{k,j}^\frac{1}{2} + D_{k,f}^\frac{1}{2} + D_{k,\rho}^\frac{1}{2}$.

\subsection{Higher-order mixed estimates on microscopic equations} 
\label{sub:mixed_estimates}


We will perform in this subsection higher-order mixed derivative estimates for both microscopic equations. In fact, we consider here the higher-order mixed derivatives on $\nm{\nabla_x^k \nabla_q^l (f_A,\,f_B)}_{\lxq}$, in which $k+l=s,\, l \ge 1$. Besides, we denote $\nabla_l^k = \nabla_x^k \nabla_q^l$ for notational simplicity.

Applying the mixed derivative operator $\nabla_l^k$ on the first and second equation of \eqref{sys:perturbative}, we get respectively,
  \begin{multline}\label{Eq:deri-xq-A}
    \p_t \nabla_l^k f_A + \nabla_l^k (u \cdot \nabla_x f_A) + \nabla_l^k (\nabla_x u q \nabla_q f_A) = \nabla_l^k \mathcal{L}_A f_A \\ + \nabla_l^k [\nabla_x u q \nabla_q U_A(\sqrt{M_A} + \tfrac{1}{2} f_A)] + \nabla_l^k r_A ,
  \end{multline}
and
  \begin{multline}\label{Eq:deri-xq-B}
    \p_t \nabla_l^k f_B + \nabla_l^k (u \cdot \nabla_x f_B) + \nabla_l^k (\nabla_x u q \nabla_q f_B) = \nabla_l^k \mathcal{L}_B f_B \\ + \nabla_l^k [\nabla_x u q \nabla_q U_B (\sqrt{M_B} + \tfrac{1}{2} f_B)] + \nabla_l^k r_B .
  \end{multline}

It follows that,
  \begin{align}
    \skpa{\p_t \nabla_l^k f_A}{\nabla_l^k f_A} = \frac{1}{2} \frac{\d}{\d t} \nm{\nabla_l^k f_A}_{\lxq}^2,
  \end{align}
and by a discussion by cases, that
  \begin{align}
    & \skpa{\nabla_l^k (u \cdot \nabla_x f_A)}{\nabla_l^k f_A} \\\no
    & = \skpa{[\nabla_x^k,\, u \cdot \nabla_x] \nabla_q^l f_A}{q \nabla_l^k f_A} + \skpa{u \cdot \nabla_x ( \nabla_l^k f_A)}{\nabla_l^k f_A}, \\\no
    & = \sum_{\substack{k_1+k_2=k\\ 1 \le k_1 \le k \le s-1,\,}} \skpa{\nabla_x^{k_1} u \nabla_x \nabla_l^{k_2} f_A}{\nabla_l^k f_A} \\\no
    & = \skpa{\nabla_x^{k} u \nabla_x \nabla_q^l f_A}{\nabla_l^k f_A} + \sum_{\substack{k_1+k_2=k\\1 \le k_1 \le k-1 \le s-2}} \skpa{\nabla_x^{k_1} u \nabla_l^{k_2+1} f_A}{\nabla_l^k f_A} \\\no
    & \le \Big(\abs{\nabla_x^{k} u}_{L^4_x} \nm{\nabla_x \nabla_q^l f_A}_{L^4_x L^2_q}
            + \sum_{\substack{k_1+k_2=k\\1 \le k_1 \le k-1 \le s-2}} \abs{\nabla_x^{k_1} u}_{L^\infty_x} \nm{\nabla_l^{k_2+1} f_A}_{L^2_x L^2_q}
          \Big)
          \nm{\nabla_l^k f_A}_{L^2_x L^2_q} \\\no
    & \ls \abs{u}_{\hx{s}} \nm{\nabla_q^l f_A}_{\hxq{k}} \nm{\nabla_l^k f_A}_{\lxq}.
  \end{align}

As for the third term in the left-hand side of \eqref{Eq:deri-xq-A}, we notice that
  \begin{align}
    \skpa{\nabla_l^k (\nabla_x u q \nabla_q f_A)}{\nabla_l^k f_A}
    = \skpa{\nabla_x^k (\nabla_x u q \nabla_q^l f_A)}{\nabla_q \nabla_l^k f_A}
      + \skpa{\nabla_x^k (\nabla_x u \nabla_q^l f_A)}{\nabla_l^k f_A}.
  \end{align}
A similar discussion by cases as above ensures that
  \begin{align}
    & \skpa{\nabla_x^k (\nabla_x u q \nabla_q^l f_A)}{\nabla_q \nabla_l^k f_A} \\\no
    & = \skpa{\nabla_x^{k+1} u q \nabla_q^l f_A}{\nabla_q \nabla_l^k f_A}
        + \skpa{\nabla_x^{k} u q \nabla_l^1 f_A}{\nabla_q \nabla_l^k f_A}
        + \sum_{\substack{k_1+k_2=k\\k_1 \le k-2}} \skpa{\nabla_x^{k_1+1} u q \nabla_l^{k_2} f_A}{\nabla_q \nabla_l^k f_A} \\\no
    & \le \Big(\abs{\nabla_x^{k+1} u}_{L^2_x} \nm{q \nabla_q^l f_A}_{L^\infty_x L^2_q}
            + \abs{\nabla_x^{k} u}_{L^4_x} \nm{q \nabla_l^1 f_A}_{L^4_x L^2_q} \\[4pt] \no
          & \hspace*{5.2cm} + \sum_{\substack{k_1+k_2=k\\ k_1 \le k-2}} \abs{\nabla_x^{k_1+1} u}_{L^\infty_x} \nm{q \nabla_l^{k_2} f_A}_{L^2_x L^2_q} \Big)
          \nm{\nabla_q \nabla_l^k f_A}_{L^2_x L^2_q} \\[4pt] \no
    & \ls \abs{u}_{\hx{s}} \nm{q \nabla_q^l f_A}_{\hxq{k}} \nm{\nabla_q \nabla_l^k f_A}_{\lxq},
  \end{align}
and
  \begin{align}
    & \skpa{\nabla_x^k (\nabla_x u \nabla_q^l f_A)}{\nabla_l^k f_A} \\\no
    & = \skpa{\nabla_x^{k+1} u \nabla_q^l f_A}{\nabla_l^k f_A}
        + \skpa{\nabla_x^{k} u \nabla_l^1 f_A}{\nabla_l^k f_A}
        + \sum_{\substack{k_1+k_2=k\\k_1 \le k-2}} \skpa{\nabla_x^{k_1+1} u \nabla_l^{k_2} f_A}{\nabla_l^k f_A} \\\no
    & \ls \abs{u}_{\hx{s}} \nm{\nabla_q^l f_A}_{\hxq{k}} \nm{\nabla_l^k f_A}_{\lxq},
  \end{align}
then we get
  \begin{align}
    & \skpa{\nabla_l^k (\nabla_x u q \nabla_q f_A)}{\nabla_l^k f_A} \\\no
    & \ls \abs{u}_{\hx{s}} \nm{q \nabla_q^l f_A}_{\hxq{k}} \nm{\nabla_q \nabla_l^k f_A}_{\lxq}
          + \abs{u}_{\hx{s}} \nm{\nabla_q^l f_A}_{\hxq{k}} \nm{\nabla_l^k f_A}_{\lxq}.
  \end{align}

We turn to consider the right-hand side terms in \eqref{Eq:deri-xq-A}. Firstly, the expression of Fokker-Planck operator $\mathcal{L}_\alpha f_\alpha = \Delta_q f_\alpha + \tfrac{1}{2} \Delta_q U_\alpha f_\alpha - \tfrac{1}{4} |\nabla_q U_\alpha|^2 f_\alpha$ yields the commutator relation:
  \begin{align}
    \nabla_l^k \L_A f_A = \L_A \nabla_l^k f_A + \sum_{\substack{l_1+l_2 = l\\ l_1 \ge 1}} \tfrac{1}{2} \nabla_q^{l_1} (\Delta_q U_A - \tfrac{1}{2} |\nabla_q U_A|^2) \nabla_q^{l_2} \nabla_x^k f_A,
  \end{align}
Recall the assumptions on potentials \eqref{asmp-2}, $\abs{\nabla_q^{l_1} (\Delta_q U_A - \tfrac{1}{2} |\nabla_q U_A|^2)} \ls |\nabla_q U_A|^2 + 1$, then this relation implies
  \begin{align}
    & \skpa{\nabla_l^k \L_A f_A}{\nabla_l^k f_A} \\\no
    & = \skpa{\L_A \nabla_l^k f_A}{\nabla_l^k f_A} + \skpa{[\nabla_l^k, \L_A] f_A}{\nabla_l^k f_A} \\\no
    & \ge - \nm{\nabla_q(\tfrac{\nabla_l^k f_A}{\sqrt{M_A}})}_M^2 \\\no
      & \qquad \quad + C \sum_{0\le l_2 \le l-1} \Big( \nm{\nabla_q U_A \nabla_{l_2}^k f_A}_{\lxq} \nm{\nabla_q U_A \nabla_l^k f_A}_{\lxq} + \nm{\nabla_{l_2}^k f_A}_{\lxq} \nm{\nabla_l^k f_A}_{\lxq} \Big).
  \end{align}

As for the penultimate term in right-hand side of \eqref{Eq:deri-xq-A}, the potential assumption \eqref{asmp-2} $\abs{\nabla_q^{l_1} (q \nabla_q U_A)} \ls \abs{q \nabla_q U_A}$ enables us to get
  \begin{align}
    & \skpa{\nabla_l^k (\nabla_x u q \nabla_q U_A f_A)}{\nabla_l^k f_A} \\\no
    & = \sum_{\substack{k_1+k_2=k\\l_1+l_2 =l,\, l \ge 1}} \skpa{\nabla_x^{k_1+1} u\, \nabla_q^{l_1} (q \nabla_q U_A)\, \nabla_{l_2}^{k_2} f_A}{\nabla_l^k f_A} \\\no
    & \ls \sum_{l_2 \le l} \Big(\abs{\nabla_x^{k+1} u}_{L^2_x} \nm{q \nabla_q^{l_2} f_A}_{L^\infty_x L^2_q}
            + \abs{\nabla_x^{k} u}_{L^4_x} \nm{q \nabla_{l_2}^1 f_A}_{L^4_x L^2_q} \\[4pt] \no
          & \hspace*{5.2cm} + \sum_{\substack{k_1+k_2=k\\ k_1 \le k-2}} \abs{\nabla_x^{k_1+1} u}_{L^\infty_x} \nm{q \nabla_{l_2}^{k_2} f_A}_{L^2_x L^2_q} \Big)
          \nm{\nabla_q U_A \nabla_l^k f_A}_{L^2_x L^2_q} \\[4pt] \no
    & \ls \abs{u}_{\hx{s}} \sum_{l_2 \le l} \nm{q \nabla_q^{l_2} f_A}_{\hxq{k}} \nm{\nabla_q U_A \nabla_l^k f_A}_{\lxq},
  \end{align}
and the potential assumption \eqref{asmp-2} $\nm{\nabla_q^l (q \nabla_q U_A \sqrt{M_A})}_{L^2_q} \le C$ ensures that
  \begin{align}
    \skpa{\nabla_l^k (\nabla_x u q \nabla_q U_A \sqrt{M_A})}{\nabla_l^k f_A}
    & \le \abs{\nabla_x^{k+1} u}_{L^2_x} \nm{\nabla_q^l (q \nabla_q U_A \sqrt{M_A})}_{L^2_q} \nm{\nabla_l^k f_A}_{\lxq} \\\no
    & \ls \abs{\nabla_x u}_{\hx{k}} \nm{\nabla_l^k f_A}_{\lxq}.
  \end{align}

Now we are in a position to deal with the reaction contributions, which can be reformulated as
  \begin{align}
    & \skpa{\nabla_l^k r_A}{\nabla_l^k f_A} + \skpa{\nabla_l^k r_B}{\nabla_l^k f_B} \\\no
    & = - \skpa{\nabla_l^k (f_A - 2 f_B \sqrt{M_B} - f_B^2)}{\nabla_l^k f_A} + 2 \skpa{\nabla_l^k (f_A \sqrt{M_B} - 2 f_B {M_B} - f_B^2 \sqrt{M_B})}{\nabla_l^k f_A} \\\no
    & = - \nm{\nabla_l^k f_A - 2 \nabla_l^k f_B \sqrt{M_B}}_{\lxq}^2
          + \iint \nabla_l^k(f_B^2) (\nabla_l^k f_A - 2 \nabla_l^k f_B \sqrt{M_B}) \d q \d x
          + \sum_{l_2 \ge 1} \sum_{i=1}^4 \texttt{Re}_i,
  \end{align}
where the remainder terms $\texttt{Re}_i$'s can be controlled by
  \begin{align}
    \abs{\texttt{Re}_1} & = \abs{2 \iint \nabla_{l_1}^k f_B \nabla_q^{l_2} \sqrt{M_B}\, \nabla_l^k f_A \d q \d x} \\\no
    & \le 2 \nm{\nabla_{l_1}^k f_B}_{\lxq} \nm{\nabla_q^{l_2} \sqrt{M_B}}_{L^\infty_q} \nm{\nabla_l^k f_A}_{\lxq} \\\no
    & \ls \nm{\nabla_{l_1}^k f_B}_{\lxq} \nm{\nabla_l^k f_A}_{\lxq},
  \end{align}
and similarly,
  \begin{align*}
    \abs{\texttt{Re}_2} & = \abs{2 \iint \nabla_{l_1}^k f_A \nabla_q^{l_2} \sqrt{M_B}\, \nabla_l^k f_B \d q \d x}
      \ls \nm{\nabla_{l_1}^k f_A}_{\lxq} \nm{\nabla_l^k f_B}_{\lxq}, \\\no
    \abs{\texttt{Re}_3} & = \abs{-4 \iint \nabla_{l_1}^k f_B \nabla_q^{l_2} {M_B}\, \nabla_l^k f_B \d q \d x}
      \ls \nm{\nabla_{l_1}^k f_B}_{\lxq} \nm{\nabla_l^k f_B}_{\lxq}, \\\no
    \abs{\texttt{Re}_4} & = \abs{-2 \iint \nabla_{l_1}^k (f_B^2) \nabla_q^{l_2} \sqrt{M_B}\, \nabla_l^k f_B \d q \d x} \ls \nm{\nabla_{l_1}^k (f_B^2)}_{\lxq} \nm{\nabla_l^k f_B}_{\lxq}.
  \end{align*}
We then get
  \begin{align}
    & \skpa{\nabla_l^k r_A}{\nabla_l^k f_A} + \skpa{\nabla_l^k r_B}{\nabla_l^k f_B} \\\no
    & \ls - \nm{\nabla_l^k f_A - 2 \nabla_l^k f_B \sqrt{M_B}}_{\lxq}^2
          + \nm{\nabla_l^k (f_B^2)}_{\lxq} \nm{\nabla_l^k f_A - 2 \nabla_l^k f_B \sqrt{M_B}}_{\lxq} \\\no
        & \quad  + \sum_{l_1 \le l-1} \Big(\nm{\nabla_{l_1}^k f_B}_{\lxq} \nm{\nabla_l^k f_A}_{\lxq}
            + \nm{\nabla_{l_1}^k f_A}_{\lxq} \nm{\nabla_l^k f_B}_{\lxq} \\\no
          & \hspace*{2.1cm} + \nm{\nabla_{l_1}^k f_B}_{\lxq} \nm{\nabla_l^k f_B}_{\lxq}
            + \nm{\nabla_{l_1}^k (f_B^2)}_{\lxq} \nm{\nabla_l^k f_B}_{\lxq} \Big)
      \\\no
    & \ls - \nm{\nabla_l^k f_A - 2 \nabla_l^k f_B \sqrt{M_B}}_{\lxq}^2
          + \nm{f_B}_{H^2_x H^2_q} \nm{f_B}_{\hxqt{s}} \nm{\nabla_l^k f_A - 2 \nabla_l^k f_B \sqrt{M_B}}_{\lxq}
        \\\no & \quad
          + \Big[\nm{f_B}_{\hxqt{s-1}} \nm{\nabla_l^k f_A}_{\lxq}
          + \big( \nm{f_A}_{\hxqt{s-1}} + \nm{f_B}_{\hxqt{s-1}} + \nm{f_B}_{H^2_x H^2_q} \nm{f_B}_{\hxqt{s-1}} \big) \nm{\nabla_l^k f_B}_{\lxq} \Big].
  \end{align}

Therefore, we have derived the higher-order mixed estimates on microscopic equations, for index satisfying $k+l=s,\, l \ge 1$, that
  \begin{align}\label{esm:mix-xq-AB}
    & \frac{1}{2} \frac{\d}{\d t} \nm{\nabla_l^k (f_A,\, f_B)}_{\lxq}^2
          + \nm{ \nabla_q \left( \tfrac{\nabla_l^k f_A}{\sqrt{M_A}} \right) }_M^2
          + \nm{ \nabla_q \left( \tfrac{\nabla_l^k f_B}{\sqrt{M_B}} \right) }_M^2
          + \nm{\nabla_l^k f_A - 2 \nabla_l^k f_B \sqrt{M_B}}_{\lxq}^2 \no\\
    & \ls |u|_{\hx{s}} \nm{\nabla_q^l f_A}_{\hxq{k}} \nm{\nabla_l^k f_A}_{\lxq}
        + |u|_{\hx{s}} \nm{q \nabla_q^l f_A}_{\hxq{k}} \nm{\nabla_q \nabla_l^k f_A}_{\lxq} \\\no
    & \quad + \sum_{0\le l_2 \le l-1} \left[
                  \nm{\nabla_q U_A \nabla_q^{l_2} f_A}_{\hxq{k}} \nm{\nabla_q U_A \nabla_l^k f_A}_{\lxq}
                + \nm{\nabla_q^{l_2} f_A}_{\hxq{k}} \nm{\nabla_l^k f_A}_{\lxq} \right] \\\no
    & \quad + |\nabla_x u|_{\hx{k}} \nm{\nabla_l^k f_A}_{\lxq}
            + |u|_{\hx{s}} \sum_{0\le l_2 \le l} \nm{q \nabla_q^{l_2} f_A}_{\hxq{k}} \nm{\nabla_q U_A \nabla_l^k f_A}_{\lxq}
      \\\no
    & \quad + |u|_{\hx{s}} \nm{\nabla_q^l f_B}_{\hxq{k}} \nm{\nabla_l^k f_B}_{\lxq}
            + |u|_{\hx{s}} \nm{q \nabla_q^l f_B}_{\hxq{k}} \nm{\nabla_q \nabla_l^k f_B}_{\lxq} \\\no
    & \quad + \sum_{0\le l_2 \le l-1} \left[
                  \nm{\nabla_q U_B \nabla_q^{l_2} f_B}_{\hxq{k}} \nm{\nabla_q U_B \nabla_l^k f_B}_{\lxq}
                + \nm{\nabla_q^{l_2} f_B}_{\hxq{k}} \nm{\nabla_l^k f_B}_{\lxq} \right] \\\no
    & \quad + |\nabla_x u|_{\hx{k}} \nm{\nabla_l^k f_B}_{\lxq}
            + |u|_{\hx{s}} \sum_{0\le l_2 \le l} \nm{q \nabla_q^{l_2} f_B}_{\hxq{k}} \nm{\nabla_q U_B \nabla_l^k f_B}_{\lxq}
    \\\no
    & \quad + k_1 \nm{f_B}_{H^2_x H^2_q} \nm{f_B}_{\hxqt{s}} \nm{\nabla_l^k f_A - 2 \nabla_l^k f_B \sqrt{M_B}}_{\lxq} \\\no
    & \quad + \sum_{l_1 \le l-1} \left[ \nm{\nabla_{l_1}^k f_B}_{\hxq{k}} \nm{\nabla_l^k f_A}_{\lxq}
      + \nm{\nabla_{l_1}^k f_A}_{\hxq{k}} \nm{\nabla_l^k f_B}_{\lxq}
      + \nm{\nabla_{l_1}^k f_B}_{\hxq{k}} \nm{\nabla_l^k f_B}_{\lxq} \right]  \\\no
    & \quad + \nm{f_B}_{H^2_x H^2_q} \nm{f_B}_{\hxqt{s}} \nm{\nabla_l^k f_B}_{\lxq}.
  \end{align}
This yields that
  \begin{align}
    & \frac{\d}{\d t} E_{s,\mix} + D_{s,\mix} \\\no
    & \ls \eta E_{s,u}^\frac{1}{2} (D_{s-1,\mix} + D_{s-1,f} + D_{s-1,\rho}) \\\no
    & \quad + \eta^\frac{1}{2} E_{s,u}^\frac{1}{2} (D_{s,\mix}^\frac{1}{2} + \eta^\frac{1}{2} D_{s-1,f}^\frac{1}{2} + \eta^\frac{1}{2} D_{s-1,\rho}^\frac{1}{2}) (D_{s-1,\mix}^\frac{1}{2} + D_{s-1,f}^\frac{1}{2} + D_{s-1,\rho}^\frac{1}{2}) \\\no
    & \quad + \eta^\frac{1}{2} (D_{s-1,f}^\frac{1}{2} + D_{s-1,\rho}^\frac{1}{2}) (D_{s-1,\mix}^\frac{1}{2} + \eta^\frac{1}{2} D_{s-1,f}^\frac{1}{2} + \eta^\frac{1}{2} D_{s-1,\rho}^\frac{1}{2})
      \\\no
    & \quad + \sum_{\substack{1 \le l_2 \le l-1\\l\ge 2}}
                (\eta^\frac{l-l_2}{2} D_{s-1,\mix}^\frac{1}{2} + \eta^\frac{l}{2} D_{s-1,f}^\frac{1}{2} + \eta^\frac{l}{2} D_{s-1,\rho}^\frac{1}{2})
                (D_{s,\mix}^\frac{1}{2} + \eta^\frac{l}{2} D_{s-1,f}^\frac{1}{2} + \eta^\frac{l}{2} D_{s-1,\rho}^\frac{1}{2}) \\\no
    & \quad + \sum_{\substack{1 \le l_2 \le l-1\\l\ge 2}}
            (\eta^\frac{l-l_2+1}{2} D_{s-2,\mix}^\frac{1}{2} + \eta^\frac{l}{2} D_{s-1,f}^\frac{1}{2} + \eta^\frac{l}{2} D_{s-1,\rho}^\frac{1}{2})
            (\eta^\frac{1}{2} D_{s-1,\mix}^\frac{1}{2} + \eta^\frac{l}{2} D_{s-1,f}^\frac{1}{2} + \eta^\frac{l}{2} D_{s-1,\rho}^\frac{1}{2}) \\\no
    & \quad + \sum_{\substack{k+l =s\\l\ge 1}} \eta^\frac{l+1}{2} D_{k,u}^\frac{1}{2} (D_{s-1,\mix}^\frac{1}{2} + D_{s-1,f}^\frac{1}{2} + D_{s-1,\rho}^\frac{1}{2}) \\\no
    & \quad + \sum_{0\le l_2 \le l} E_{s,u}^\frac{1}{2}
                    (\eta^\frac{l-l_2}{2} D_{s,\mix}^\frac{1}{2} + \eta^\frac{l}{2} D_{s-1,f}^\frac{1}{2} + \eta^\frac{l}{2} D_{s-1,\rho}^\frac{1}{2})
                    (D_{s,\mix}^\frac{1}{2} + \eta^\frac{l}{2} D_{s-1,f}^\frac{1}{2} + \eta^\frac{l}{2} D_{s-1,\rho}^\frac{1}{2}) \\\no
    & \quad + k_1 \sum_{\substack{k+l =s\\l\ge 1}} \eta^\frac{l}{2} \nm{f_B}_{H^2_x H^2_q} \nm{f_B}_{\hxqt{s}} D_{s,\mix,r}^\frac{1}{2} \\\no
    & \quad + \sum_{\substack{k+l =s\\l\ge 1}} \eta^\frac{l+1}{2} \nm{f_B}_{H^2_x H^2_q} \nm{f_B}_{\hxqt{s}}
            (D_{s-1,\mix}^\frac{1}{2} + \eta^\frac{l-1}{2} D_{s-1,f}^\frac{1}{2} + \eta^\frac{l-1}{2} D_{s-1,\rho}^\frac{1}{2}) .
  \end{align}

\subsection{Additional dissipation on number densities} 
\label{sub:dissipation_density}


In the above process of energy estimates, the weighted Poincar\'e inequalities are frequently used to help us control the energy terms by dissipation terms, in which the mean value functions of number densities $\rho_{\alpha}(t,\,x)$ are involved. This causes a new difficulty in closing the estimates because these lower order terms seem to be pure energy terms without some natural dissipation effects.

As we described before, the mean value terms are due to the non-conservation of each distribution function $f_{\alpha}$, or speaking specifically, due to the reaction terms. This suggests us to focus on the reaction contributions by extracting them from the micro-macro structure of kinetic equations.

Indeed, the micro-macro coupling will vanish by performing integrations over configuration space due to their nice divergence structure of Smoluchowski equation, as a correspondence to the conservation law part in the left-hand side of the non-conservative kinematics \eqref{eq:kinematic}. More precisely, we perform an integration on both two kinetic equations in the perturbative system \eqref{sys:perturbative} with respect to the weighted measure $\sqrt{M_\alpha} \d q$, then get equations of perturbative number densities $\rho_{\alpha}(t,\,x) = \int f_{\alpha} \sqrt{M_\alpha} \d q = \agl{g_{\alpha}}_M$, involving (almost) pure reaction contributions, as follows,
  \begin{align}\label{eq:number-density}
  \begin{cases}
    \p_t \rho_A + u \cdot \nabla_x \rho_A = - \agl{r}, \\
    \p_t \rho_B + u \cdot \nabla_x \rho_B = 2 \agl{r}, \\
  \end{cases}
  \end{align}
where
  \begin{align}\label{eq:integ-react}
    \agl{r(q)}(t,\,x) = \agl{\sqrt{M_A} (f_A - 2 \sqrt{M_B} f_B - f_B^2)}.
  \end{align}

We point out there holds the conservation relation for the total perturbations,
  \begin{align}
    2 \rho_A + \rho_B = 0.
  \end{align}
which is actually a direct consequence of \eqref{eq:conserv-total} revealing in a macroscopic formulation.

Now we consider the higher-order derivatives estimates for number density $\abs{\nabla_x^s \rho_A}_{\lx}$ (or equivalently, $\abs{\nabla_x^s \rho_B}_{\lx}$). Applying higher-order derivative operator $\nabla_x^s$ on the first equation of \eqref{eq:number-density}, and taking $L^2_x$-inner product with $\nabla_x^s \rho_A$, we can get
  \begin{align}\label{eq:estimate-rho}
    \frac{1}{2} \frac{\d}{\d t} |\nabla_x^s \rho_A|^2_{\lx}
    = - \skpa{\nabla_x^s(u \cdot \nabla_x \rho_A)}{\nabla_x^s \rho_A} - \skpa{\nabla_x^s r}{\nabla_x^s \rho_A}.
  \end{align}

Recall the definition of projection on the kernel of linear Fokker-Planck operator $\P_0 g = \rho^g$,
the linear part in above reaction contribution \eqref{eq:integ-react} can be recast as,
  \begin{align}
    \agl{\sqrt{M_A} (f_A - 2 \sqrt{M_B} f_B)}
    & = \agl{g_A M_A - 2 \sqrt{M_A M_B} \cdot (g_B \sqrt{M_B})} \\\no
    & = \rho_A - 2 \agl{(\P_0 g_B + \P_0^\perp g_B) M_B^2} \\\no
    & = (\rho_A - 2 \agl{M_A} \rho_B) + \agl{(\P_0^\perp g_B) M_B^2} \\\no
    & = (1 + 4 \agl{M_A}) \rho_A + \agl{(\P_0^\perp g_B) M_B^2}.
  \end{align}
due to the above relation $2 \rho_A + \rho_B = 0$. In fact, the expression will lead to an additional dissipation term, $(1 + 4 \agl{M_A}) \abs{\nabla_x^s \rho_A}_{\lx}^2$, in energy estimate for number density \eqref{eq:estimate-rho}.

Notice that we have established the partial coercivity estimates on the kernel orthogonal part of Fokker-Planck operator in \eqref{eq:partial-coer}, or in other words, the Poincar\'e inequality without a mean value function, which enables us to get
  \begin{align}
    \iint (\P_0^\perp g_B) M_B^2 \cdot \rho_A \d q \d x
    & \le \agl{M_B^3}^\frac{1}{2} \left( \iint \abs{\P_0^\perp g_B}^2 M_B \d q \d x \right)^\frac{1}{2} \abs{\rho_A}_{\lx} \\\no
    & \le \lambda_0^{-1} \agl{M_B^3}^\frac{1}{2} \skpa{-\A g_B}{g_B}_M^\frac{1}{2} \cdot \abs{\rho_A}_{\lx} \\\no
    & \ls \skpa{-\L f_B}{f_B}^\frac{1}{2} \cdot \abs{\rho_A}_{\lx},
  \end{align}
and similarly, for higher-order derivatives case,
  \begin{align}
    \iint (\P_0^\perp \nabla_x^s g_B) M_B^2 \cdot \nabla_x^s \rho_A \d q \d x
      \ls \skpa{-\L \nabla_x^s f_B}{\nabla_x^s f_B}^\frac{1}{2} \cdot \abs{\nabla_x^s \rho_A}_{\lx}.
  \end{align}

Since we have
  \begin{align}
    \skpa{\nabla_x^s(u \cdot \nabla_x \rho_A)}{\nabla_x^s \rho_A}
    & = \skpa{[\nabla_x^s, u \cdot \nabla_x] \rho_A)}{\nabla_x^s \rho_A} + \skpa{u \cdot \nabla_x (\nabla_x^s \rho_A)}{\nabla_x^s \rho_A} \\\no
    & \ls |u|_{\hx{s}} |\rho_A|_{\hx{s}} |\nabla_x^s \rho_A|_{\lx}.
  \end{align}
we can derive finally the higher-order estimates for number density, that
  \begin{align}\label{eq:dissipation-density}
    & \frac{1}{2} \frac{\d}{\d t} |\nabla_x^s \rho_A|^2_{\lx} + (4 \agl{M_A} + 1) |\nabla_x^s \rho_A|^2_{\lx} \\\no
    & \ls |u|_{\hx{s}} |\rho_A|_{\hx{s}} |\nabla_x^s \rho_A|_{\lx}
    + |\nabla_x^s \rho_A|_{\lx} \nm{\nabla_q ( \tfrac{\nabla_x^s f_B}{\sqrt{M_B}})}_M 
    + |\nabla_x^s \rho_A|_{\lx} \nm{f_B}_{H^2_x H^2_q} \nm{f_B}_{\hxq{s}}.
  \end{align}
This can also be written as
  \begin{align}
    \frac{\d}{\d t} E_{s,\rho} + D_{s,\rho}
    \ls E_{s,u}^\frac{1}{2} \D_{s,\rho} + \D_{s,\rho}^\frac{1}{2} \D_{s,f}^\frac{1}{2}
        + \nm{f_B}_{H^2_x H^2_q} \nm{f_B}_{\hxq{s}} \D_{s,\rho}^\frac{1}{2} .
  \end{align}

\subsection{Closing the \emph{a-priori} estimate} 
\label{sub:closing_estimates}


Before combining all the above estimates to obtain the closed \emph{a-priori} estimate, we are required firstly to control product terms appeared in the form $\nm{f_B}_{H^2_x H^2_q} \nm{f_B}_{X}$ with $X$ denoting another norm.

\smallskip\noindent\underline{\textbf{Necessary bound on quadratic terms}}:

We consider the following products involving $\nm{f_B}_{H^2_x H^2_q}$, which is viewed as a dissipation part. More precisely, for $k\le 2$,
  \begin{align}
    \nm{\nabla_x^k f_B}_{L^2_x H^2_q}^2
    & = \nm{\nabla_q^2 \nabla_x^k f_B}_{\lxq}^2 + \nm{\nabla_q \nabla_x^k f_B}_{\lxq}^2 + \nm{\nabla_x^k f_B}_{\lxq}^2 \\\no
    & \ls \nm{ \nabla_q \big( \tfrac{\nabla_x^k \nabla_q f_B}{\sqrt{M_B}} \big) }_M^2
          + \nm{ \nabla_q \big( \tfrac{\nabla_x^k f_B}{\sqrt{M_B}} \big) }_M^2
          + \abs{\nabla_x^k \rho_B}_{\lx}^2,
  \end{align}
as a consequence, it holds,
  \begin{align}
    \eta \nm{f_B}_{H^2_x H^2_q}^2 \ls D_{3,\mix} + \eta D_{2,f} + \eta D_{2,\rho}.
  \end{align}
This enables us to get
  \begin{align}
   \eta^\frac{1}{2} \nm{f_B}_{H^2_x H^2_q} \nm{f_B}_{\hxq{s}}
   \ls E_{s,f}^\frac{1}{2} (D_{3,\mix}^\frac{1}{2} + \eta^\frac{1}{2} D_{2,f}^\frac{1}{2} + \eta^\frac{1}{2} D_{2,\rho}^\frac{1}{2}).
  \end{align}

Similar process implies that
  \begin{align}
    \eta \nm{qf_B}_{H^2_x H^2_q} \nm{qf_B}_{\hxq{s'}}
      & \ls E_{s',j}^\frac{1}{2} (D_{4,\mix}^\frac{1}{2} + \eta D_{2,f}^\frac{1}{2} + \eta D_{2,\rho}^\frac{1}{2}),
   \\
    \eta^\frac{s+1}{2} \nm{f_B}_{H^2_x H^2_q} \nm{f_B}_{\hxqt{s}}
      & \ls E_{4,\mix}^\frac{1}{2} (D_{s-1,\mix}^\frac{1}{2} + \eta^\frac{s-1}{2} D_{s,f}^\frac{1}{2} + \eta^\frac{s-1}{2} D_{s,\rho}^\frac{1}{2}).
  \end{align}

\smallskip\noindent\underline{\textbf{Closing the estimates}}:

We are now ready to close the desired estimates. Recall the definitions of energy and dissipation functionals introduced in Proposition \ref{prop:a-priori}
  \begin{align}
    \E_s^\eta & = E_{s} + \eta E_{s,\rho} + \eta E_{s,\mix} + \eta^2 E_{s',j} , \\
    \D_s^\eta & = D_{s} + \eta D_{s,\rho} + \eta D_{s,\mix} + \eta^2 D_{s',j} ,
  \end{align}
then we can infer, respectively, for pure spatial estimates, that
  \begin{align}
    \frac{\d}{\d t} E_{s} + D_{s} 
      \ls \eta^{-1} (E_{s,u}^\frac{1}{2} + E_{s,f}^\frac{1}{2}) \D_s^\eta + \eta^{-\frac{3}{2}} [(\eta^2 E_{s',j})^\frac{1}{2} + (\eta E_{s,\mix})^\frac{1}{2}] \, \D_s^\eta.
  \end{align}
for pure spatial estimates on first order moment, that
  \begin{align}
    \frac{\d}{\d t} (\eta^2 E_{s',j}) + (\eta^2 D_{s',j})
      \ls \eta \D_s^\eta + \eta D_{s',u}^\frac{1}{2} (\D_s^\eta)^\frac{1}{2}
          + E_{s,u}^\frac{1}{2} \D_s^\eta + \eta^{-\frac{1}{2}} (\eta^2 E_{s',j})^\frac{1}{2}\, \D_s^\eta.
  \end{align}
for mixed spatial-configurational estimates, that
  \begin{align}
    \frac{\d}{\d t} (\eta E_{s,\mix}) + (\eta D_{s,\mix})
    \ls \eta^\frac{1}{2} \D_s^\eta
          + & \sum_{\substack{k+l =s\\l\ge 1}} \eta^\frac{l}{2} D_{k,u}^\frac{1}{2} ( \D_s^\eta )^\frac{1}{2}
          \\\no &
          + E_{s,u}^\frac{1}{2} \D_s^\eta + \eta^{-\frac{s+1}{2}} (\eta E_{s,\mix})^\frac{1}{2}\, \D_s^\eta.
  \end{align}
and for additional dissipation estimates on number density, that
  \begin{align}
    \frac{\d}{\d t} E_{s,\rho}^\eta + D_{s,\rho}^\eta
    \ls \eta^\frac{1}{2} \D_s^\eta + E_{s,u}^\frac{1}{2} \D_s^\eta + \eta^{-\frac{1}{2}} (E_{s,f}^\eta)^\frac{1}{2}\, \D_s^\eta.
  \end{align}

Finally, we can derive by combining all the above estimates together, that
  \begin{align}\label{esm:a-priori-var}
     \frac{\d}{\d t} \E_s^\eta + \D_s^\eta
     \ls \eta^\frac{1}{2} \D_s^\eta + \eta^{-\frac{s}{2}} (\E_s^\eta)^\frac{1}{2}\, \D_s^\eta
          + \eta D_{s',u}^\frac{1}{2} (\D_s^\eta)^\frac{1}{2} + \sum_{\substack{k+l =s\\l\ge 1}} \eta^\frac{l}{2} D_{k,u}^\frac{1}{2} (\D_s^\eta)^\frac{1}{2},
  \end{align}
which yields immediately the desired \emph{a-priori} estimate \eqref{esm:a-priori} that,
  \begin{align}
    \frac{\d}{\d t} \E_s^\eta + \D_s^\eta
    \ls ( \eta^\frac{1}{2} + \eta^{-\frac{s}{2}} \eps^\frac{1}{2} ) \D_s^\eta,
  \end{align}
provided $\E_s^\eta \le \eps$. This completes the whole proof of Proposition \ref{prop:a-priori}. \qed

\begin{remark}
Note that, we will use in the following texts the formulation \eqref{esm:a-priori-var} for the {a-priori} estimate, in which the last two terms in the right-hand side are kept. This is very important to help justify the global-in-time existence by a continuity argument, as we will see below.
\end{remark}

\section{Completion of the Proof for Global Existence} 
\label{sec:completion_of_proof}


To prove Theorem \ref{thm:main}, a local existence result is required. We state it through the following lemma, 
the proof of which relies on a standard iterating scheme and compactness method and is thus omitted here. The readers are referred to, e.g. \cite{LZZ-04cpde,Ren-91sima,JLZ-18sima}.

\begin{lemma}[Local existence] \label{lemm:local-exis}
 Let $s\ge 5$ and $\epsilon_0>0$ be a small positive constant. Assume $\E(0) \le \epsilon/2$ for any constant $\epsilon \in (0,\, \epsilon_0)$, then there exists a time $T_*>0$ such that the perturbative system \eqref{sys:perturbative} of two-species micro-macro model for wormlike micellar solutions admits a unique local classical solution $(u,\, f_A,\ f_B) \in L^\infty(0, T_*; H^s_x \times H^s_{x,q} \times H^s_{x,q})$, satisfying,
  \begin{align}
    \E(t) \le \epsilon.
  \end{align}
Moreover, the positivity $\Psi_{\alpha} = M_\alpha + \sqrt{M_\alpha} f_{\alpha} >0$ is preserved on $[0,\, T_*]$ if it is valid for initial data $\Psi_{\alpha,0} = M_\alpha + \sqrt{M_\alpha} f_{\alpha,0} >0$.
\end{lemma}

We now in a position of proving our main result Theorem \ref{thm:main}. By the equivalence of energy functionals $\E_s \sim \E_s^\eta$, it follows from \eqref{eq:initial-small} that, $\E_s^\eta(0) \le 2 \E_s(0) \le 2\eps$. Lemma \ref{lemm:local-exis} then ensures that, there exists a local solution with the maximal existing time $T_*>0$, moreover, for $t \in [0,\,T_*]$,
  \begin{align}\label{eq:local-contradict}
    \E_s^\eta(t) \le 4 \eps.
  \end{align}
We claim that $T_* = + \infty$. For this, we prove by contradiction.

Recall the \emph{a-priori} estimate \eqref{esm:a-priori-var} obtained in Proposition \ref{prop:a-priori}, then there exists some positive constant $C>0$, such that the following inequality holds
  \begin{align}\label{eq:a-priori-var-1}
    \frac{\d}{\d t} \E_s^\eta + \D_s^\eta
     \le C (\eta^\frac{1}{2} + \eta^{-\frac{s}{2}} \eps^\frac{1}{2}) \D_s^\eta
          + C \eta \abs{\nabla_x u}_{\hx{s'}} (\D_s^\eta)^\frac{1}{2} + C \sum_{\substack{k+l =s\\l\ge 1}} \eta^\frac{l}{2} \abs{\nabla_x u}_{\hx{k}} (\D_s^\eta)^\frac{1}{2},
  \end{align}
for any $t \in [0,\,T_*]$.

As for the last two terms in above expression, the H\"older inequality yields that
  \begin{align}
    \eta \abs{\nabla_x u}_{\hx{s'}} (\D_s^\eta)^\frac{1}{2} + \sum_{\substack{k+l =s\\l\ge 1}} \eta^\frac{l}{2} \abs{\nabla_x u}_{\hx{k}} (\D_s^\eta)^\frac{1}{2}
    \le C \eta^\frac{1}{2} \D_s^\eta + C \eta^\frac{1}{2} \abs{\nabla_x u}_{L^2_x}^2,
  \end{align}
where we have used the interpolation inequality of Sobolev spaces, for $k\le s'=s-1$ and some small parameter $\delta>0$, that
  \begin{align*}
    \abs{\nabla_x u}_{\hx{k}}^2 \le \delta \abs{\nabla_x u}_{\hx{s}}^2 + C_{\delta} \abs{\nabla_x u}_{L^2_x}^2,
  \end{align*}
and the simple fact $\abs{\nabla_x u}_{\hx{s}}^2 \le \D_s^\eta$.

Therefore, by virtue of the small assumptions on $\eps$ and $\eta$, it follows from \eqref{eq:a-priori-var-1} that,
  \begin{align}
    \frac{\d}{\d t} \E_s^\eta + \tfrac{1}{2} \D_s^\eta \le C \eta^\frac{1}{2} \abs{\nabla_x u}_{L^2_x}^2,
  \end{align}
which implies immediately, for $t \in [0,\, T_*]$,
  \begin{align}\label{eq:continuity-argu}
    \E_s^\eta(t) + \int_0^t \tfrac{1}{2} \D_s^\eta(\tau) \d \tau
    \le \E_s^\eta(0) + C \eta^\frac{1}{2} \int_0^t \abs{\nabla_x u}_{L^2_x}^2 \d \tau.
  \end{align}

We now turn back to the basic energy-dissipation law \eqref{law:energy-dissipation} by employing the notion of relative entropy. We introduce the (relative) energy density functions, for $\alpha=A,B$,
  \begin{align}\label{eq:relative density}
    E_{\alpha}(\Psi_{\alpha} | M_{\alpha}) = \Psi_{\alpha} \ln \frac{\Psi_{\alpha}}{M_{\alpha}} - \Psi_{\alpha} + M_{\alpha}.
  \end{align}
It is easy to check its non-negativity that $E_{\alpha}(\Psi_{\alpha} | M_{\alpha}) = M_{\alpha} h(\tfrac{\Psi_{\alpha} - M_{\alpha}}{M_{\alpha}}) \ge 0$, by considering a convex function $h(z) = (1+z)\ln (1+z) - z$ which is defined on $(-1,+\infty)$.

As a consequence, replacing in \eqref{law:energy-dissipation} the original energy density $\widetilde E_{\alpha} = \Psi_{\alpha} (\ln \Psi_{\alpha} + U_{\alpha} -1)$ by the relative energy density $E_{\alpha}(\Psi_{\alpha} | M_{\alpha})$, we get the following refined energy-dissipation law for the two-species micro-macro model \eqref{sys:2sp-mic-mac}:
  \begin{align}\label{law:refined-basic}
    \frac{\d}{\d t} & \left\{ \int_\Omega \tfrac{1}{2}|u|^2 \d x
      + \lambda \sum_{\alpha=A,B} \iint_{\Omega \times \mathbb{R}^3} E_{\alpha}(\Psi_{\alpha} | M_{\alpha}) \d q \d x \right\} \\\no
    =\ & - \int_\Omega \mu |\nabla_x u |^2 \d x
         - \lambda \sum_{\alpha=A,B} \iint_{\Omega \times \mathbb{R}^3} \Psi_\alpha \abs{\nabla_q (\ln \Psi_\alpha + U_\alpha)}^2 \d q \d x
         \\\no
       & - \lambda \iint_{\Omega \times \mathbb{R}^3} (k_1 \Psi_A - k_2 \Psi_B^2) (\ln \tfrac{\Psi_A}{\Psi_B^2} + U_A - 2 U_B) \d q \d x,
  \end{align}
which immediately yields
  \begin{align}
    \frac{\d}{\d t} \left\{ \int_\Omega \tfrac{1}{2}|u|^2 \d x
      + \lambda \sum_{\alpha=A,B} \iint_{\Omega \times \mathbb{R}^3} E_{\alpha}(\Psi_{\alpha} | M_{\alpha}) \d q \d x \right\}
    \le 0.
  \end{align}
By integrating over $[0,\,t]$, this inequality implies
  \begin{align}
    \tfrac{1}{2}|u|_{L^2_x}^2(t) & + \int_0^t \mu \abs{\nabla_x u}_{L^2_x}^2(\tau) \d \tau \\\no
    & \le \tfrac{1}{2}|u_0|_{L^2_x}^2
    + \lambda \iint_{\Omega \times \mathbb{R}^3} [E_A(\Psi_{A,0}| M_A)  + E_B(\Psi_{B,0}| M_B)] \d q \d x
    \le \eps,
  \end{align}
due to the small assumption on initial data \eqref{eq:initial-small-entropy}. Then \eqref{eq:continuity-argu} reduces to
  \begin{align*}\label{eq:continuity-argu2}
    \E_s^\eta(t) + \int_0^t \tfrac{1}{2} \D_s^\eta(\tau) \d \tau
    \le 2\eps + C \eta^\frac{1}{2} \eps
    \le 3 \eps,
  \end{align*}
which is valid for all $t \in [0,\,T_*]$ and especially for $t= T_*$. This contradicts \eqref{eq:local-contradict} unless it holds that $T_* = +\infty$.

Recalling the equivalence between energy functionals with and without the fixed parameter $\eta$, the above energy inequality in turn yields the global energy bound \eqref{eq:glob-energy-bdd}. This completes the whole proof of Theorem \ref{thm:main}. \qed



\section*{Acknowledgments}

This work is partially supported by the National Science Foundation (USA) grants NSF DMS-1950868, the United States-Israel Binational Science Foundation (BSF) \#2024246 (C. Liu, Y. Wang), and a grant from the National Natural Science Foundation of China No. 11871203 (T.-F. Zhang). This work was done when T.-F. Zhang visited the Department of Applied Mathematics at Illinois Institute of Technology during 2019-2020, he would like to acknowledge the hospitality of IIT and the sponsorship of the China Scholarship Council, under the State Scholarship Fund No. 201906415023. The authors would like to thank Prof. Ning Jiang for helpful discussions.


\bibliographystyle{aomvar}
\bibliography{wlm-refs}

\providecommand{\bysame}{\leavevmode\hbox to3em{\hrulefill}\thinspace}
\providecommand{\noopsort}[1]{}
\providecommand{\mr}[1]{\href{http://www.ams.org/mathscinet-getitem?mr=#1}{MR~#1}}
\providecommand{\zbl}[1]{\href{http://www.zentralblatt-math.org/zmath/en/search/?q=an:#1}{Zbl~#1}}
\providecommand{\jfm}[1]{\href{http://www.emis.de/cgi-bin/JFM-item?#1}{JFM~#1}}
\providecommand{\arxiv}[1]{\href{http://www.arxiv.org/abs/#1}{arXiv~#1}}
\providecommand{\doi}[1]{\url{https://doi.org/#1}}
\providecommand{\MR}{\relax\ifhmode\unskip\space\fi MR }
\providecommand{\MRhref}[2]{%
  \href{http://www.ams.org/mathscinet-getitem?mr=#1}{#2}
}
\providecommand{\href}[2]{#2}
\begin{thebibliography}{10}

\bibitem{BCAH-87b2}
\bgroup{}R.~B. Bird\egroup{}, \bgroup{}C.~F. Curtiss\egroup{}, \bgroup{}R.~C.
  Armstrong\egroup{}, and \bgroup{}O.~Hassager\egroup{}, \emph{Dynamics of
  polymeric liquids, Volume 2: Kinetic theory}, second ed., John Wiley \& Sons,
  New York, 1987.

\bibitem{Cat87}
\bgroup{}M.~E. Cates\egroup{}, Reptation of living polymers: dynamics of
  entangled polymers in the presence of reversible chain-scission reactions,
  \emph{Macromolecules} \textbf{20} (1987) no.~9, 2289--2296.

\bibitem{CM-01sima}
\bgroup{}J.-Y. Chemin\egroup{} and \bgroup{}N.~Masmoudi\egroup{}, About
  lifespan of regular solutions of equations related to viscoelastic fluids,
  \emph{SIAM J. Math. Anal.} \textbf{33} (2001) no.~1, 84--112.

\bibitem{DL-97arma}
\bgroup{}P.~Degond\egroup{} and \bgroup{}M.~Lemou\egroup{}, Dispersion
  relations for the linearized {F}okker-{P}lanck equation,  \emph{Arch. Ration.
  Mech. Anal.} \textbf{138} (1997) no.~2, 137--167.

\bibitem{DLP-17jns}
\bgroup{}P.~Degond\egroup{}, \bgroup{}J.-G. Liu\egroup{}, and \bgroup{}R.~L.
  Pego\egroup{}, Coagulation-fragmentation model for animal group-size
  statistics,  \emph{J. Nonlinear Sci.} \textbf{27} (2017) no.~2, 379--424.

\bibitem{DE-86b}
\bgroup{}M.~Doi\egroup{} and \bgroup{}S.~F. Edwards\egroup{}, \emph{The theory
  of polymer dynamics}, \textbf{73}, Oxford University Press, Oxford, UK, 1986.

\bibitem{DFT-10cmp}
\bgroup{}R.~Duan\egroup{}, \bgroup{}M.~Fornasier\egroup{}, and
  \bgroup{}G.~Toscani\egroup{}, A kinetic flocking model with diffusion,
  \emph{Comm. Math. Phys.} \textbf{300} (2010) no.~1, 95--145.

\bibitem{ELZ-04cmp}
\bgroup{}W.~E\egroup{}, \bgroup{}T.~Li\egroup{}, and
  \bgroup{}P.~Zhang\egroup{}, Well-posedness for the dumbbell model of
  polymeric fluids,  \emph{Comm. Math. Phys.} \textbf{248} (2004) no.~2,
  409--427.

\bibitem{GCB-13thermo}
\bgroup{}N.~Germann\egroup{}, \bgroup{}L.~Cook\egroup{}, and
  \bgroup{}A.~Beris\egroup{}, Nonequilibrium thermodynamic modeling of the
  structure and rheology of concentrated wormlike micellar solutions,  \emph{J.
  Non-Newton. Fluid Mech.} \textbf{196} (2013), 51--57.

\bibitem{GKL-18notes}
\bgroup{}M.-H. Giga\egroup{}, \bgroup{}A.~Kirshtein\egroup{}, and
  \bgroup{}C.~Liu\egroup{}, Variational modeling and complex fluids,  in
  \emph{Handbook of mathematical analysis in mechanics of viscous fluids},
  Springer, Cham, 2018, pp.~73--113.

\bibitem{GST-13siap}
\bgroup{}T.~Goudon\egroup{}, \bgroup{}M.~Sy\egroup{}, and \bgroup{}L.~M.
  Tin\'{e}\egroup{}, A fluid-kinetic model for particulate flows with
  coagulation and breakup: stationary solutions, stability, and hydrodynamic
  regimes,  \emph{SIAM J. Appl. Math.} \textbf{73} (2013) no.~1, 401--421.

\bibitem{Guo-02cpam}
\bgroup{}Y.~Guo\egroup{}, The {V}lasov-{P}oisson-{B}oltzmann system near
  {M}axwellians,  \emph{Comm. Pure Appl. Math.} \textbf{55} (2002) no.~9,
  1104--1135.

\bibitem{Guo-03invent}
\bgroup{}Y.~Guo\egroup{}, The {V}lasov-{M}axwell-{B}oltzmann system near
  {M}axwellians,  \emph{Invent. Math.} \textbf{153} (2003) no.~3, 593--630.

\bibitem{Guo-06cpam}
\bgroup{}Y.~Guo\egroup{}, Boltzmann diffusive limit beyond the
  {N}avier-{S}tokes approximation,  \emph{Comm. Pure Appl. Math.} \textbf{59}
  (2006) no.~5, 626--687.

\bibitem{Hu-18jde}
\bgroup{}X.~Hu\egroup{}, Global existence of weak solutions to two dimensional
  compressible viscoelastic flows,  \emph{J. Differential Equations}
  \textbf{265} (2018) no.~7, 3130--3167.

\bibitem{HL-16cpam}
\bgroup{}X.~Hu\egroup{} and \bgroup{}F.-H. Lin\egroup{}, Global solutions of
  two-dimensional incompressible viscoelastic flows with discontinuous initial
  data,  \emph{Comm. Pure Appl. Math.} \textbf{69} (2016) no.~2, 372--404.

\bibitem{HLL-18notes}
\bgroup{}X.~Hu\egroup{}, \bgroup{}F.-H. Lin\egroup{}, and
  \bgroup{}C.~Liu\egroup{}, Equations for viscoelastic fluids,  in
  \emph{Handbook of mathematical analysis in mechanics of viscous fluids},
  Springer, Cham, 2018, pp.~1045--1073.

\bibitem{HCDL-08krm}
\bgroup{}Y.~Hyon\egroup{}, \bgroup{}J.~A. Carrillo\egroup{},
  \bgroup{}Q.~Du\egroup{}, and \bgroup{}C.~Liu\egroup{}, A maximum entropy
  principle based closure method for macro-micro models of polymeric materials,
   \emph{Kinet. Relat. Models} \textbf{1} (2008) no.~2, 171--184.

\bibitem{JLZ-18sima}
\bgroup{}N.~Jiang\egroup{}, \bgroup{}Y.~Liu\egroup{}, and \bgroup{}T.-F.
  Zhang\egroup{}, Global classical solutions to a compressible model for
  micro-macro polymeric fluids near equilibrium,  \emph{SIAM J. Math. Anal.}
  \textbf{50} (2018) no.~4, 4149--4179.

\bibitem{JLL-04jfa}
\bgroup{}B.~Jourdain\egroup{}, \bgroup{}T.~Leli\`evre\egroup{}, and
  \bgroup{}C.~Le~Bris\egroup{}, Existence of solution for a micro-macro model
  of polymeric fluid: the {FENE} model,  \emph{J. Funct. Anal.} \textbf{209}
  (2004) no.~1, 162--193.

\bibitem{KLLM-20m2an}
\bgroup{}P.~Knopf\egroup{}, \bgroup{}K.~F. Lam\egroup{},
  \bgroup{}C.~Liu\egroup{}, and \bgroup{}S.~Metzger\egroup{}, Phase-field
  dynamics with transfer of materials: {T}he {C}ahn--{H}illard equation with
  reaction rate dependent dynamic boundary conditions,  \emph{ESAIM Math.
  Model. Numer. Anal.} \textbf{55} (2021) no.~1, 229--282.

\bibitem{LbL-12scm}
\bgroup{}C.~Le~Bris\egroup{} and \bgroup{}T.~Leli\`evre\egroup{}, Micro-macro
  models for viscoelastic fluids: modelling, mathematics and numerics,
  \emph{Sci. China Math.} \textbf{55} (2012) no.~2, 353--384.

\bibitem{LLZ-08arma}
\bgroup{}Z.~Lei\egroup{}, \bgroup{}C.~Liu\egroup{}, and
  \bgroup{}Y.~Zhou\egroup{}, Global solutions for incompressible viscoelastic
  fluids,  \emph{Arch. Ration. Mech. Anal.} \textbf{188} (2008) no.~3,
  371--398.

\bibitem{LZZ-04cpde}
\bgroup{}T.~Li\egroup{}, \bgroup{}H.~Zhang\egroup{}, and
  \bgroup{}P.~Zhang\egroup{}, Local existence for the dumbbell model of
  polymeric fluids,  \emph{Comm. Partial Differential Equations} \textbf{29}
  (2004) no.~5-6, 903--923.

\bibitem{LZ-07cms}
\bgroup{}T.~Li\egroup{} and \bgroup{}P.~Zhang\egroup{}, Mathematical analysis
  of multi-scale models of complex fluids,  \emph{Commun. Math. Sci.}
  \textbf{5} (2007) no.~1, 1--51.

\bibitem{Lin-12cpam}
\bgroup{}F.-H. Lin\egroup{}, Some analytical issues for elastic complex fluids,
   \emph{Comm. Pure Appl. Math.} \textbf{65} (2012) no.~7, 893--919.

\bibitem{LLZ-07cpam}
\bgroup{}F.-H. Lin\egroup{}, \bgroup{}C.~Liu\egroup{}, and
  \bgroup{}P.~Zhang\egroup{}, On a micro-macro model for polymeric fluids near
  equilibrium,  \emph{Comm. Pure Appl. Math.} \textbf{60} (2007) no.~6,
  838--866.

\bibitem{LZZ-08cmp}
\bgroup{}F.-H. Lin\egroup{}, \bgroup{}P.~Zhang\egroup{}, and
  \bgroup{}Z.~Zhang\egroup{}, On the global existence of smooth solution to the
  2-{D} {FENE} dumbbell model,  \emph{Comm. Math. Phys.} \textbf{277} (2008)
  no.~2, 531--553.

\bibitem{LM-00cam}
\bgroup{}P.~L. Lions\egroup{} and \bgroup{}N.~Masmoudi\egroup{}, Global
  solutions for some {O}ldroyd models of non-{N}ewtonian flows,  \emph{Chinese
  Ann. Math. Ser. B} \textbf{21} (2000) no.~2, 131--146.

\bibitem{Lc-09notes}
\bgroup{}C.~Liu\egroup{}, An introduction of elastic complex fluids: an
  energetic variational approach,  in \emph{Multi-scale phenomena in complex
  fluids}, \emph{Ser. Contemp. Appl. Math. CAM} \textbf{12}, World Sci.
  Publishing, Singapore, 2009, pp.~286--337.

\bibitem{LW-19arma}
\bgroup{}C.~Liu\egroup{} and \bgroup{}H.~Wu\egroup{}, An energetic variational
  approach for the {C}ahn-{H}illiard equation with dynamic boundary condition:
  model derivation and mathematical analysis,  \emph{Arch. Ration. Mech. Anal.}
  \textbf{233} (2019) no.~1, 167--247.

\bibitem{LYY-04physd}
\bgroup{}T.-P. Liu\egroup{}, \bgroup{}T.~Yang\egroup{}, and \bgroup{}S.-H.
  Yu\egroup{}, Energy method for {B}oltzmann equation,  \emph{Phys. D}
  \textbf{188} (2004) no.~3-4, 178--192.

\bibitem{LY-04cmp}
\bgroup{}T.-P. Liu\egroup{} and \bgroup{}S.-H. Yu\egroup{}, Boltzmann equation:
  micro-macro decompositions and positivity of shock profiles,  \emph{Comm.
  Math. Phys.} \textbf{246} (2004) no.~1, 133--179.

\bibitem{Maj-84b}
\bgroup{}A.~J. Majda\egroup{}, \emph{Compressible fluid flow and systems of
  conservation laws in several space variables}, \emph{Applied Mathematical
  Sciences} \textbf{53}, Springer-Verlag, New York, 1984.

\bibitem{Mas-08cpam}
\bgroup{}N.~Masmoudi\egroup{}, Well-posedness for the {FENE} dumbbell model of
  polymeric flows,  \emph{Comm. Pure Appl. Math.} \textbf{61} (2008) no.~12,
  1685--1714.

\bibitem{Mas-13invent}
\bgroup{}N.~Masmoudi\egroup{}, Global existence of weak solutions to the {FENE}
  dumbbell model of polymeric flows,  \emph{Invent. Math.} \textbf{191} (2013)
  no.~2, 427--500.

\bibitem{Mas-18notes}
\bgroup{}N.~Masmoudi\egroup{}, Equations for polymeric materials,  in
  \emph{Handbook of mathematical analysis in mechanics of viscous fluids},
  Springer, Cham, 2018, pp.~973--1005.

\bibitem{Old1950}
\bgroup{}J.~G. Oldroyd\egroup{}, On the formulation of rheological equations of
  state,  \emph{Proc. Roy. Soc. London Ser. A} \textbf{200} (1950), 523--541.

\bibitem{Ons1931-1}
\bgroup{}L.~Onsager\egroup{}, Reciprocal relations in irreversible processes.
  {I}.,  \emph{Phys. Rev.} \textbf{37} (1931) no.~4, 405--426.

\bibitem{Ons1931-2}
\bgroup{}L.~Onsager\egroup{}, Reciprocal relations in irreversible processes.
  {II}.,  \emph{Phys. Rev.} \textbf{38} (1931) no.~12, 2265--2279.

\bibitem{OP-74arma}
\bgroup{}G.~F. Oster\egroup{} and \bgroup{}A.~S. Perelson\egroup{}, Chemical
  reaction dynamics. {I}. {G}eometrical structure,  \emph{Arch. Ration. Mech.
  Anal.} \textbf{55} (1974), 230--274.

\bibitem{Ren-91sima}
\bgroup{}M.~Renardy\egroup{}, An existence theorem for model equations
  resulting from kinetic theories of polymer solutions,  \emph{SIAM J. Math.
  Anal.} \textbf{22} (1991) no.~2, 313--327.

\bibitem{Ray1871}
\bgroup{}H.~J.~W. Strutt (L.~Rayleigh)\egroup{}, Some general theorems relating
  to vibrations,  \emph{Proc. Lond. Math. Soc.} \textbf{4} (1871/73), 357--368.

\bibitem{Tay-PDE3}
\bgroup{}M.~E. Taylor\egroup{}, \emph{Partial differential equations {III}.
  {N}onlinear equations}, second ed., \emph{Applied Mathematical Sciences}
  \textbf{117}, Springer, New York, 2011.

\bibitem{VCM07}
\bgroup{}P.~A. Vasquez\egroup{}, \bgroup{}G.~H. McKinley\egroup{}, and
  \bgroup{}L.~P. Cook\egroup{}, A network scission model for wormlike micellar
  solutions: I. model formulation and viscometric flow predictions,  \emph{J.
  Non-Newton. Fluid Mech.} \textbf{144} (2007) no.~2, 122--139.

\bibitem{WG86}
\bgroup{}P.~Waage\egroup{} and \bgroup{}C.~M. Gulberg\egroup{}, Studies
  concerning affinity,  \emph{J. Chem. Educ.} \textbf{63} (1986) no.~12, 1044.

\bibitem{WLLE-20pre}
\bgroup{}Y.~Wang\egroup{}, \bgroup{}C.~Liu\egroup{}, \bgroup{}P.~Liu\egroup{},
  and \bgroup{}B.~Eisenberg\egroup{}, Field theory of reaction-diffusion: {L}aw
  of mass action with an energetic variational approach,  \emph{Phys. Rev. E}
  \textbf{102} (2020) no.~6, 062147.

\bibitem{LWZ-21nnfm}
\bgroup{}Y.~Wang\egroup{}, \bgroup{}T.-F. Zhang\egroup{}, and
  \bgroup{}C.~Liu\egroup{}, A two species micro-macro model of wormlike
  micellar solutions and its maximum entropy closure approximations: {A}n
  energetic variational approach,  \emph{J. Non-Newton. Fluid Mech.}
  \textbf{293} (2021), 104559.

\bibitem{WXL-13arma}
\bgroup{}H.~Wu\egroup{}, \bgroup{}X.~Xu\egroup{}, and \bgroup{}C.~Liu\egroup{},
  On the general {E}ricksen-{L}eslie system: {P}arodi's relation,
  well-posedness and stability,  \emph{Arch. Ration. Mech. Anal.} \textbf{208}
  (2013) no.~1, 59--107.

\bibitem{YJZ-10ccp}
\bgroup{}H.~Yu\egroup{}, \bgroup{}G.~Ji\egroup{}, and
  \bgroup{}P.~Zhang\egroup{}, A nonhomogeneous kinetic model of liquid crystal
  polymers and its thermodynamic closure approximation,  \emph{Commun. Comput.
  Phys.} \textbf{7} (2010) no.~2, 383--402.

\bibitem{ZZ-06arma}
\bgroup{}H.~Zhang\egroup{} and \bgroup{}P.~Zhang\egroup{}, Local existence for
  the {FENE}-dumbbell model of polymeric fluids,  \emph{Arch. Ration. Mech.
  Anal.} \textbf{181} (2006) no.~2, 373--400.

\end{thebibliography}

%

\end{document}